\documentclass[12pt]{amsart}
\usepackage{amsfonts}
\usepackage{amssymb}
\usepackage{xypic}

\newtheorem{theorem}{Theorem}[subsection]
\newtheorem{lem}[theorem]{Lemma}
\newtheorem{cor}[theorem]{Corollary}
\newtheorem{prop}[theorem]{Proposition}
\newtheorem{ass}[theorem]{Assertion}
\newtheorem{df}[theorem]{Definition-Proposition}

\theoremstyle{definition}
\newtheorem{definition}[theorem]{Definition}

\theoremstyle{remark}
\newtheorem{remark}[theorem]{Remark}

\theoremstyle{definition}

\theoremstyle{remark}

\theoremstyle{plain}

\def\limind{\mathop{\oalign{lim\cr
\hidewidth$\longrightarrow$\hidewidth\cr}}}
\def\boxit#1#2{\setbox1=\hbox{\kern#1{#2}\kern#1}%
\dimen1=\ht1 \advance\dimen1 by #1
\dimen2=\dp1 \advance\dimen2 by #1
\setbox1=\hbox{\vrule height\dimen1 depth\dimen2\box1\vrule}%
\setbox1=\vbox{\hrule\box1\hrule}%
\advance\dimen1 by .4pt \ht1=\dimen1
\advance\dimen2 by .4pt \dp1=\dimen2 \box1\relax}

\def\AA{{\mathbf A}}
\def\BB{{\mathbf B}}
\def\CC{{\mathbf C}}

\def\FF{{\mathbf F}}
\def\GG{{\mathbf G}}

\def\LL{{\mathbf L}}

\def\NN{{\mathbf N}}

\def\PP{{\mathbf P}}
\def\QQ{{\mathbf Q}}

\def\ZZ{{\mathbf Z}}

\def\cD{{\mathcal D}}
\def\cE{{\mathcal E}}
\def\cF{{\mathcal F}}

\def\cI{{\mathcal I}}
\def\cJ{{\mathcal J}}
\def\cK{{\mathcal K}}
\def\cL{{\mathcal L}}
\def\cM{{\mathcal M}}

\def\cO{{\mathcal O}}

\def\cV{{\mathcal V}}

\def\cZ{{\mathcal Z}}

\mathchardef\alphag="7C0B
\mathchardef\betag="7C0C
\mathchardef\gammag="7C0D
\mathchardef\deltag="7C0E
\mathchardef\varepsilong="7C22
\mathchardef\varphig="7C27
\mathchardef\psig="7C20
\mathchardef\zetag="7C10
\mathchardef\epsilong="7C0F
\mathchardef\rhog="7C1A
\mathchardef\taug="7C1C
\mathchardef\upsilong="7C1D
\mathchardef\iotag="7C13
\mathchardef\thetag="7C12
\mathchardef\pig="7C19
\mathchardef\sigmag="7C1B
\mathchardef\etag="7C11
\mathchardef\omegag="7C21
\mathchardef\kappag="7C14
\mathchardef\lambdag="7C15
\mathchardef\mug="7C16
\mathchardef\xig="7C18
\mathchardef\chig="7C1F
\mathchardef\nug="7C17
\mathchardef\varthetag="7C23
\mathchardef\varpig="7C24
\mathchardef\varrhog="7C25
\mathchardef\varsigmag="7C26
\mathchardef\Omegag="7C0A
\mathchardef\Thetag="7C02
\mathchardef\Sigmag="7C06
\mathchardef\Deltag="7C01
\mathchardef\Phig="7C08
\mathchardef\Gammag="7C00
\mathchardef\Psig="7C09
\mathchardef\Lambdag="7C03
\mathchardef\Xig="7C04
\mathchardef\Pig="7C05
\mathchardef\Upsilong="7C07

\begin{document}

\title{Motivic Igusa zeta functions}

\author{Jan Denef}
\address{University of Leuven, Department of Mathematics,
Celestijnenlaan 200B, 3001 Leuven, Belgium}
\email{Jan.Denef@wis.kuleuven.ac.be}

\author{Fran\c cois Loeser}

\address{Centre de Math\'ematiques,
Ecole Polytechnique,
F-91128 Palaiseau
(URA 169 du CNRS), {\rm and}
Institut de Math\'{e}matiques,
Universit\'{e} P. et M. Curie, Case 82,
4 place Jussieu,
F-75252 Paris Cedex 05
(UMR 9994 du CNRS)}
\email{loeser@math.polytechnique.fr}

\date{revised september 1997}

\dedicatory{To appear in Journal of Algebraic Geometry}



\maketitle

\setcounter{section}{-1}
\renewcommand{\theequation}{\thesection.\arabic{equation}}
\section{Introduction}Let $p$
be a prime number and let $K$ be a finite extension of $\QQ_p$.
Let $R$ be the valuation ring of $K$, $P$ the maximal ideal of $R$, and
$\bar K = R / P$ the residue field of $K$. Let $q$ denote the cardinality of
$\bar K$, so $\bar K \simeq \FF_q$. For $z$ in $K$,  let
${\rm ord} \, z$ denote the valuation of $z$, and  set
$\vert z \vert = q^{- {\rm ord} \, z}$.
Let $f$ be a non constant element of $K [x_1, \ldots, x_m]$.
The $p$-adic Igusa local zeta function $Z (s)$ associated to $f$
(relative to 
the trivial multiplicative character) is defined as the 
$p$-adic integral
\begin{equation}
Z (s) = \int_{R^m} \vert f (x) \vert^s \vert dx \vert,
\end{equation}
for $s \in \CC$, ${\rm Re} (s) > 0$, 
where $\vert dx \vert$ denotes the Haar measure on
$K^m$ normalized in such of way that $R^m$ is of volume 1.
For $n$ in $\NN$, set
$
Z_{n} = \{ x \in R^{m} \bigm \vert {\rm ord} \, f (x) = n\}.
$
We may express $Z (s)$ as a series
\begin{equation}
Z (s) = \sum_{n \geq 0} \, {\rm vol} \, (Z_{n}) \,
q^{-ns}.
\end{equation}
Now, if we denote by $X_{n}$ the image of $Z_{n}$ in  
$(R / P^{n + 1})^{m}$, we may rewrite the series as
\begin{equation}
Z (s) = \sum_{n \geq 0} \, {\rm card} \, (X_{n}) \,
q^{-ns - (n + 1) m}
\end{equation}
since ${\rm vol} \, (Z_{n}) =  {\rm card} \, (X_{n}) \, q^{- (n+ 1) m}$.

\bigskip

Let $k$ be a field of characteristic zero. M.Kontsevich recently
introduced
the concept of motivic integration \cite{K} (see also \cite{D-L2} and 2.5),
which is a $k [[t]]$-analogue of usual
$p$-adic integration. This motivic integration takes values into a
certain completion of a localisation of the Grothendieck ring
$K_{0} ({\rm Sch}_{k})$
of algebraic varieties over $k$, {\it i.e.}
reduced separated schemes of finite type over $k$,
(see 2.5 for more details). 
The ring $K_{0} ({\rm Sch}_{k})$
is
generated by symbols $[S]$, for $S$ an algebraic variety over
$k$, 
with the relations
$[S] = [S']$ if $S$ is isomorphic to $S'$,
$[S] = [S \setminus S'] + [S']$ if $S'$ is closed in $S$
and
$[S \times S'] = [S] \, [S']$.
We set $\LL := [\AA^1_k]$.
In this setting, the analogue of $R^{m}$, {\it resp.} of
$(R / P^{n + 1})^{m}$,
is the $k$-scheme $\cL (\AA^{m}_{k})$, {\it resp.} $\cL_{n} (\AA^{m}_{k})$,
which represents the functor $R \mapsto R [[t]]^{m}$, {\it resp.} the functor
$R \mapsto (R [[t]] / t^{n+1} R [[t]])^{m}$,
on
the category of $k$-algebras, {\it i.e.} the
$k$-scheme parametrizing $m$-tuples of series,
{\it resp.} of series modulo $t^{n + 1}$.
Let $f$ be a non constant element of $K [x_1, \ldots, x_m]$.
Define $Z_{n}$ as the subscheme
of $\cL
(\AA^{m}_{k})$ of series $\varphi$ such that ${\rm ord}_{t} f (\varphi) = n$
and $X_{n}$ as the image of $Z_{n}$ in 
$\cL_{n} (\AA^{m}_{k})$, viewed as a reduced subscheme. A natural
analogue of the right-hand side of (0.3), which is a series in 
$\ZZ [p^{-1}]  [[p^{-s}]]$, 
is the following series in
$K_{0} ({\rm Sch}_{k}) [\LL^{-1}]  [[\LL^{-s}]]$
\begin{equation}
Z_{{\rm geom}} (s) = \sum_{n \geq 0} \, [X_{n}] \,
\LL^{-ns - (n + 1) m}.
\end{equation}
When 
$s$ takes a fixed value in $\NN$, this series can be interpreted as a Kontsevich
integral (see 2.5).

\bigskip

More generally, $p$-adic Igusa local zeta functions involve
multiplicative characters. Let $\pi$ be a fixed uniformizing parameter
of $R$
and set
${\rm ac} (z) = z \pi^{- {\rm ord} \, z}$ for $z$ in $K$.
For any character $\alpha : R^{\times} \rightarrow \CC^{\times}$,
one defines the $p$-adic Igusa local zeta function $Z (s, \alpha)$
as the integral
\begin{equation}
Z (s, \alpha) = \int_{R^m} \alpha ({\rm ac} (f (x))) \vert f (x) \vert^s \vert dx \vert,
\end{equation}
for $s \in \CC$, ${\rm Re} (s) > 0$ (see \cite{I}, \cite{D3}).
To extend the definition of (0.4) to the more general situation involving
characters,
it is necessary to replace varieties by motives. More generally,
let  $X$ be a smooth connected separated
scheme of finite type
over $k$ (a field of characteristic zero), let
$W$ be a reduced subscheme of $X$, and
let $f : X \rightarrow \AA^1_k$ be a morphism.
In the present paper,
we define, for $\alpha$ a multiplicative character of any finite
subgroup of $k^{\times}$, motivic Igusa functions
$\displaystyle\int_W (f^s, \alpha)$. These functions live in
a power series ring
$K_0 (\cM) [[\LL^{-s}]]$ , where $K_0 (\cM)$ is a 
Grothendieck ring of Chow motives and $\LL$ is the standard
Lefschetz motive (precise definitions are given in 1.1 and 2.1), and are
defined as series quite similar to the right-hand  side of (0.4).
These motivic Igusa functions
specialize, in the $p$-adic case with good reduction,
to the usual $p$-adic Igusa local zeta functions (see 2.4).
They also specialize to the topological zeta functions
$Z_{\rm top} (s)$ introduced by the authors
in \cite{D-L1} (see 2.3). The functions
$Z_{\rm top} (s)$ are, heuristically,  obtained as a limit as $q$ goes
to 1 of $p$-adic Igusa local zeta functions (in a more provocative way,
one can say they are defined by integrals over $W (\FF_{1})$,
the ring of Witt vectors with
coefficients in
the field
with one element).

\bigskip

The content of the paper is the following.
The motivic integrals 
$\displaystyle\int_W (f^s, \alpha)$ are defined in section 2.
Their definition uses
variants for schemes
with finite group action of
recent results of Gillet-Soul\'e \cite{G-S} and
Guill\'en-Navarro \cite{G-N}
on motivic Euler characteristics of schemes 
which are given in section 1.
In Theorem 2.2.1
we give a formula for $\displaystyle\int_W (f^s, \alpha)$ in terms of an embedded
resolution of $f$, a result which implies in particular
the rationality
of $\displaystyle\int_W (f^s, \alpha)$. We then explain the relationship
with topological zeta functions,
$p$-adic Igusa local zeta functions and motivic integration.
Section 3 is devoted to functional equations. Here
$X = W = \AA^m_k$ and  $f$ is a homogenous polynomial.
In this situation we are
able to
prove a functional equation for
motivic Igusa functions when $\alpha$ is the trivial character.
For general $\alpha$ the result depends upon a conjectural statement
on motivic Euler characteristics of quotients, but we
are able to prove it holds true if one replaces the
Grothendieck group of Chow motives
by the Grothendieck group of Voevodsky's
``triangulated
category of geometrical motives'' \cite{Vo}. These functional
equations are analogues in the present setting of
the functional equations for
$p$-adic Igusa local zeta functions proved in \cite{D-M}.
In section 4 we study
the limit for $s \rightarrow  - \infty$ of
motivic Igusa functions and investigate
its relation with
nearby cycles of $f$ at the origin. Here we are
guided by analogy with \cite{D4}, where
the limit
when $s \rightarrow - \infty$ was studied
for $p$-adic Igusa local zeta functions, and showed to be related to
the trace of some liftings of the Frobenius automorphism
acting on the cohomology of Milnor fibers.
More precisely, for
$x$ a closed point of the fiber
$f^{-1} (0)$, we give a meaning to
$$\frac{\LL^{m}}{1 - \LL}
\lim_{s \rightarrow  - \infty} \int_{\{x\}} (f^{s}, \alpha).$$
Heuristically this limit is
the ``motivic incarnation''
of $\chi_{c} (i_{x}^{\ast} R \psi_{f, \alpha})$,
where
$R \psi_{f, \alpha}$ denotes the eigenspace of nearby cycles
for the eigenvalue corresponding to
the character
$\alpha$ of the semi-simple part of the monodromy.
We prove in Theorem 4.2.1 that this holds in particular
for the $\CC$-Hodge realization. As a corollary it follows that
the whole Hodge spectrum of $f$
at $x$,
which is an important invariant of singularities
(see \cite{St1},\cite{St2}), may be deduced from the knowledge of
motivic Igusa functions.

{\small\subsection*{Acknowledgements}We would like to thank O.Villamayor
for answering to our questions on equivariant resolution and F.Morel for
interesting discussions.}

\renewcommand{\theequation}{\arabic{equation}}
\setcounter{equation}{0}

\section{Grothendieck groups of Chow motives}
\subsection{Chow motives}In this section we recall
material from \cite{Kl}, \cite{Ma}, \cite{Sc}.
We fix a base field $k$, and we denote by
$\cV_{k}$ the category of smooth and projective $k$-schemes.
For  an object $X$ in $\cV_{k}$ and an integer $d$, $\cZ^{d} (X)$
denotes the free abelian group generated by irreducible subvarieties of
$X$ of codimension $d$. We define the rational Chow group
$A^{d} (X)$ as the quotient of $\cZ^{d} (X) \otimes \QQ$
modulo rational equivalence.
For $X$ and $Y$ in $\cV_{k}$,  we denote by ${\rm Corr}^{r} (X, Y)$ the
group of correspondences of degree $r$ from $X$ to $Y$. If $X$ is purely
$d$-dimensional, ${\rm Corr}^{r} (X, Y) = A^{d + r} (X \times Y)$,
and if $X = \coprod X_{i}$, 
${\rm Corr}^{r} (X, Y) = \oplus \, {\rm Corr}^{r} (X_{i}, Y)$.
The category $\cM_{k}$ of $k$-motives may be defined as follows
(cf. \cite{Sc}).
Objects of $\cM_{k}$ are triples $(X, p, n)$ where $X$ is in $\cV_{k}$,
$p$ is an idempotent (i.e. $p^{2} = p$) in ${\rm Corr}^{0} (X, X)$, and
$n$ is an integer. If $(X, p, n)$
and $(Y, q, m)$ are motives, then
$$
{\rm Hom}_{\cM_{k}} ((X, p, n), (Y, q, m))
=
q \, {\rm Corr}^{m - n} (X, Y) \, p.
$$
Composition of morphisms is given by composition of correspondences.
The category $\cM_{k}$ is  additive, $\QQ$-linear, and pseudo-abelian.
There is a natural tensor product on $\cM_{k}$, 
defined on objects by
$$
(X, p, n) \otimes (Y, q, m) = (X \times Y, p \otimes q, n + m).
$$
We denote by $h$ the functor $h : \cV_{k}^{\circ} \rightarrow
\cM_{k}$ which sends an object $X$ to $h (X) = (X, {\rm id}, 0)$
and a morphism $f : Y \rightarrow X$ to its graph in
${\rm Corr}^{0} (X, Y)$. This functor is compatible with the
tensor product and the unit motive $1 = h ({\rm Spec} \, k)$ is the identity
for the product. We denote by ${\LL}$ the Lefschetz motive
$\LL = ({\rm Spec} \, k, {\rm id}, -1)$. There is a canonical isomorphism
$h (\PP^{1}_{k}) \simeq 1 \oplus \LL$. We denote by ${}^{\vee}$ the involution
${}^{\vee} : \cM_{k}^{\circ} \rightarrow \cM_{k}$,
defined on objects by $(X, p, n)^{\vee} =
(X, {}^{t}p, d - n)$ if $X$ is purely $d$-dimensional, and as the
transpose of correspondences on morphisms. For $X$ in
$\cV_{k}$ purely of dimension $d$, $h (X)^{\vee} = h (X) \otimes
\LL^{- d}$.

Let $E$ be a field of characteristic zero. Replacing
the Chow groups
$A^{\cdot}$ by $A^{\cdot} \otimes_{\QQ} E$, one defines similarly
the category $\cM_{k, E}$ of $k$-motives with coefficients in $E$.

\subsection{Grothendieck groups of Chow motives}
Let $K_{0} (\cM_{k})$ be the Grothendieck group of the pseudo-abelian
category
$\cM_{k}$. It is also the abelian group associated to the monoid of
isomorphism
classes of motives with respect to $\oplus$.
The tensor product on $\cM_{k}$ induces 
a natural ring structure on $K_{0} (\cM_{k})$.
Let ${\rm Sch}_{k}$ be the category of schemes which are
separated and of finite type over $k$. We suppose from now that
the
characteristic of $k$ is zero.
The following result has been
proven
by Gillet and Soul\'{e} \cite{G-S}
and also by Guill\'{e}n and Navarro Aznar \cite{G-N}.

\begin{theorem}Let $k$ be a field
of characteristic 0.
There exists a unique map
$$\chi_{c} : {\rm Ob} {\rm Sch}_{k} \longrightarrow K_{0} (\cM_{k})$$
such that
\begin{enumerate}
\item[(1)] If $X$ is smooth and projective, $\chi_{c} (X)$ is equal to
$[h (X)]$,  the image of $h (X)$ in $K_{0} (\cM_{k})$.
\item[(2)] If $Y$ is a closed subscheme in a scheme
$X$, $$\chi_{c} (X \setminus Y) = \chi_{c} (X) - \chi_{c} (Y).$$
\end{enumerate}
\end{theorem}

Let us remark that $\chi_{c} (\AA^{1}_{k}) = \LL$.
Also, for $X$ and $Y$ in ${\rm Ob} {\rm Sch}_{k}$, we have
$\chi_{c} (X \times Y) = \chi_{c} (X) \otimes \chi_{c} (Y)$.

The following result, due to Guill\'{e}n and Navarro Aznar \cite{G-N},
gives, dually, the existence of motivic Euler characteristics without supports.
A proper relative isomorphism
$(\tilde X, \tilde Y) \rightarrow (X, Y)$ consists of the following
data~:
a proper morphism $f : \tilde X \rightarrow X$ between objects of
${\rm Sch}_{k}$, reduced closed subschemes $\tilde Y$ and
$Y$ of $\tilde X$ and
$X$ respectively, such that $\tilde Y$ is the preimage of
$Y$ in $\tilde X$, and such that the restriction of $f$ to
$\tilde X \setminus \tilde Y$ is an isomorphism
onto $X \setminus Y$.

\begin{theorem}Let $k$ be a field
of characteristic 0.
There exists a map
$$\chi : {\rm Ob} \, {\rm Sch}_{k} \longrightarrow K_{0} (\cM_{k})$$
such that
\begin{enumerate}
\item[(1)] If $X$ is smooth and projective,
$\chi (X) = [h (X)]$.
\item[(2)] If $(\tilde X, \tilde Y) \rightarrow (X, Y)$ is a proper
relative isomorphism,
$$\chi (X) = \chi (\tilde X) + \chi (Y) - \chi (\tilde Y).$$
\item[(3)]If $Y$ is a smooth divisor in a smooth scheme $X$,
$$\chi (X \setminus Y) = \chi (X) - \chi (Y) \otimes \LL.$$
\item[(4)] If $X$ is a smooth scheme purely of dimension $d$,
$$\chi (X)^{\vee} = \chi_{c} (X) \otimes \LL^{-d}$$
\end{enumerate}
Furthermore, $\chi$ is determined by conditions
(1)-(3).
\end{theorem}

The Euler characteristics $\chi_{c}$ and $\chi$ are compatible with
realization
functors, in particular with Euler characteristics of mixed Hodge
structures on
cohomology with compact support and cohomology, respectively.
By additivity $\chi_{c}$ may be naturally extended to constructible sets.

\begin{remark}We expect, but do not know how to prove, that
$K_{0} (\cM_{k, E})$
has no $(\LL - 1)$-torsion. This assertion is implied by the
conjectural existence (cf. \cite{Sc} p.185) of additive functors
$h^{\leq j} : \cM_{k, E} \rightarrow \cM_{k, E}$,
$j \in \ZZ$, such that for any $X$ in 
$\cM_{k, E}$, the $h^{\leq j} (X)$ form a filtration of $X$ with
$h^{\leq -k} (X) = 0$, $h^{\leq k} (X) = X$ for some $k$, and
$h^{\leq j} (\LL X) = \LL \, h^{\leq j - 2} (X)$
for all $j$. Indeed, for $A$ in $K_{0} (\cM_{k, E})$ the 
$h^{\leq j} (A)$ are well defined in $K_{0} (\cM_{k, E})$
and the relation $(\LL - 1) A = 0 $ implies 
$h^{\leq j} (A) = \LL \, h^{\leq j - 2} (A)$, whence $A = 0$.
A similar argument also shows, without using any conjecture, that the
\'{e}tale realization and  the Hodge realization kill all
$(\LL - 1)^{i}$-torsion in $K_{0} (\cM_{k, E})$, for each $i$
in $\NN$.
\end{remark}

\subsection{Finite group action}Let $G$ be a finite abelian
group and let $\hat G$ be its complex character group.
We denote by 
$\cV_{k, G}$ the category of smooth and projective $k$-schemes
with $G$-action.
Let $E$ be a subfield of $\CC$
containing all the roots of unity of order dividing $|G|$. 
For $X$ in $\cV_{k, G}$ and $g$ in $G$, we denote by $[g]$ the
correspondence given by the graph of multiplication by $g$.

For $\alpha$ in $\hat G$ 
we consider the idempotent 
$$
f_{\alpha} := \vert G \vert^{-1} \sum_{g \in G}
\alpha^{-1} (g)
[g]
$$
in 
${\rm Corr}^{0} (X, X) \otimes E$, and
we denote by
$h (X, \alpha)$ the motive
$(X, f_{\alpha}, 0)$ in $\cM_{k, E}$.
Clearly, for $X$ in $\cV_{k, G}$ purely of dimension
$d$
and $\alpha$ in $\hat G$, we have 
$h (X, \alpha)^{\vee} = h (X, \alpha^{-1}) \otimes \LL^{-d}$.
We will denote by ${\rm Sch}_{k, G}$ the category of  
separated schemes of finite type over $k$ with $G$-action
satisfying the following condition:
the $G$-orbit of any closed point of $X$ is contained in
an affine open subscheme. This condition is clearly satisfied
for $X$ quasiprojective and insures the existence of
$X / G$  
as a scheme. Objects of ${\rm Sch}_{k, G}$ will
be called $G$-schemes.

We will need the following variants of Theorems 1.2.1 and 1.2.2. 
They are proved in the
appendix 
as a consequence of \cite{G-N} and \cite{V2}. Theorem 1.3.2 will only be used
in section 3.

\begin{theorem}Let $k$ be a field
of characteristic 0.
There exists a unique map
$$\chi_{c} : {\rm Ob} {\rm Sch}_{k, G} \times \hat G
\longrightarrow K_{0} (\cM_{k, E})$$
such that
\begin{enumerate}
\item[(1)] If $X$ is smooth and projective
with $G$-action, for any character $\alpha$,
$\chi_{c} (X, \alpha) = [h (X, \alpha)]$.
\item[(2)] If $Y$ is a closed $G$-stable
subscheme in a $G$-scheme
$X$, for any character $\alpha$,
$$\chi_{c} (X \setminus Y, \alpha) = \chi_{c} (X, \alpha) - \chi_{c} (Y, \alpha).$$
\item[(3)] If $X$ is 
a $G$-scheme, $U$ and $V$ are $G$-invariant open
subschemes of $X$, for any character $\alpha$,
$$
\chi_c (U \cup V, \alpha)
=
\chi_c (U, \alpha) +
\chi_c (V, \alpha) -
\chi_c (U \cap V, \alpha).
$$
\end{enumerate}
Furthermore, $\chi_{c}$ is determined by conditions
(1)-(2).
\end{theorem}

By a proper relative isomorphism of $G$-schemes
$(\tilde X, \tilde Y) \rightarrow (X, Y)$ we mean the following
data~:
a proper morphism $f : \tilde X \rightarrow X$ of
$G$-schemes, reduced closed $G$-stable
subschemes $\tilde Y$ and
$Y$ of $\tilde X$ and
$X$ respectively, such that $\tilde Y$ is the preimage of
$Y$ in $\tilde X$, and such that the restriction of $f$ to
$\tilde X \setminus \tilde Y$ is an $G$-isomorphism
onto $X \setminus Y$.

\begin{theorem}Let $k$ be a field
of characteristic 0.
There exists a map
$$\chi : {\rm Ob} \, {\rm Sch}_{k, G} \times \hat G
\longrightarrow K_{0} (\cM_{k, E})$$
such that
\begin{enumerate}
\item[(1)] If $X$ is smooth and projective with $G$-action,
$\chi (X, \alpha) = [h (X, \alpha)]$.
\item[(2)] If $(\tilde X, \tilde Y) \rightarrow (X, Y)$ is a proper
relative isomorphism of $G$-schemes,
$$\chi (X, \alpha) = \chi (\tilde X, \alpha) + \chi (Y, \alpha) - \chi (\tilde Y, \alpha).$$
\item[(3)]If $Y$ is a smooth $G$-invariant divisor in a smooth
$G$-scheme $X$,
$$\chi (X \setminus Y, \alpha) = \chi (X, \alpha) - \chi (Y, \alpha) \otimes \LL.$$
\item[(4)] If $X$ is a smooth $G$-scheme purely of dimension $d$,
$$\chi (X, \alpha)^{\vee} = \chi_{c} (X, \alpha^{-1}) \otimes \LL^{-d}.$$
\item[(5)]If $X$ is a proper $G$-scheme,
$$\chi_{c} (X, \alpha) = \chi (X, \alpha).$$
\item[(6)] If $X$ is 
a $G$-scheme and  $U$ and $V$ are $G$-invariant open
subschemes of $X$, then, for any character $\alpha$,
$$
\chi (U \cup V, \alpha)
=
\chi (U, \alpha) +
\chi (V, \alpha) -
\chi (U \cap V, \alpha).
$$
\end{enumerate}
Furthermore, $\chi$ is determined by conditions
(1)-(3).
\end{theorem}

\begin{prop}
Let $k$ be a field
of characteristic 0.
\begin{enumerate}
\item[(1)]For any $X$ in ${\rm Ob} {\rm Sch}_{k, G}$,
$$
\chi_{c} (X) = \sum_{\alpha \in \hat G} \chi_{c} (X, \alpha).
$$
\item[(2)]Let $X$ be in ${\rm Ob} {\rm Sch}_{k, G}$. Assume the
$G$-action factors through a quotient $G \rightarrow H$.
If $\alpha$ is not in the image of $\hat H \rightarrow \hat G$,
then $\chi_{c} (X, \alpha) = 0$.
\item[(3)]Let $X$ and $Y$ be in ${\rm Ob} {\rm Sch}_{k, G}$
and let $G$ act diagonally on $X \times Y$.
Then
$$
\chi_{c} (X \times Y, \alpha) = \sum_{\beta \in \hat G} \chi_{c} (X,
\beta) \,
\chi_{c} (Y, \alpha \beta^{-1}).
$$
\end{enumerate}
\end{prop}

\begin{proof}If $X$ is smooth and
projective with $G$-action, ${\rm Id} = \bigoplus_{\alpha \in \hat
G} f_{\alpha}$, so $h (X) = \bigoplus_{\alpha \in \hat
G} h (X, \alpha)$. It is a direct
verification that, if the $G$-action factors through a quotient
$G \rightarrow H$ and
$\alpha$ is not in the image of $\hat H \rightarrow \hat G$,
then $f_{\alpha} = 0$.
If $Y$ is another smooth and projective scheme
with $G$-action, then $$f_{\alpha} (X \times Y) = \sum_{\beta \in \hat G}
f_{\beta} (X) \otimes f_{\alpha \beta^{-1}} (Y).$$
Assertions (1), (2) and (3) follow by additivity of 
$\chi_{c} (\quad, \alpha)$.\end{proof}

\subsection{Motivic Kummer sheaves}We fix an
integer $d \geq 1$. We denote by $\mu_{d} (k)$
the group of $d$-roots of 1 in $k$ and by $\zeta_{d}$ a fixed primitive
$d$-th root of unity in $\CC$.
We assume from now on that $\mu_{d} (k)$
is of order $d$.

Let $f : X \rightarrow \GG_{m, k}$ be a morphism in ${\rm Sch}_{k}$.
For any character $\alpha$ of order $d$ of
$\mu_{d} (k)$, one may define an element $[X, f^{\ast}\cL_{\alpha}]$ of
$K_{0} (\cM_{k, \QQ[\zeta_{d}]})$ as follows.

The morphism $[d] : \GG_{m, k} \rightarrow \GG_{m, k}$
given by $x \mapsto x^{d}$ is a Galois covering with Galois group
$\mu_{d} (k)$. We consider the fiber product
\begin{equation*}\xymatrix{
\widetilde X_{f, d} \ar[d] \ar[r] & X\ar[d]^{f}\\
\GG_{m, k} \ar[r]^<<<<{[d]}& \GG_{m, k}.
}
\end{equation*}
The scheme $\widetilde X_{f, d}$ is endowed with an action of $\mu_{d} (k)$,
so we can define
$$[X, f^{\ast}\cL_{\alpha}]:= \chi_{c} (\widetilde X_{f, d}, \alpha).$$

\begin{lem}Let $f : X \rightarrow \GG_{m, k}$ and
$g : X \rightarrow \GG_{m, k}$
be morphisms in ${\rm Sch}_{k}$.
For any character $\alpha$ of order $d$ of
$\mu_{d} (k)$, the following holds
$$
[X, (f^{d} g)^{\ast} \cL_{\alpha}] =
[X, g^{\ast}\cL_{\alpha}].
$$
In particular,
$[X, f^{d\ast}\cL_{\alpha}]
= \chi_{c} (X).$
\end{lem}

\begin{proof} The morphism
$(x, t) \mapsto (x, t f^{-1} (x))$ induces
an isomorphism of $\mu_{d} (k)$-schemes
$\widetilde X_{f^{d} g, d} \simeq 
\widetilde X_{g, d}$.
For the last assertion remark that
the fiber product $\widetilde X_{1, d}$ is isomorphic 
as a $\mu_{d} (k)$-scheme to the product $X \times \mu_{d} (k)$.\end{proof}

\begin{lem}Let $f : X \rightarrow \GG_{m, k}$ and
$g : Y \rightarrow X$
be morphisms in ${\rm Sch}_{k}$. Assume that $g$ is a 
locally trivial
fibration
for the Zariski topology with fiber $Z$.
For any character $\alpha$ of order $d$ of
$\mu_{d} (k)$, the following holds
$$
[Y, (f \circ g)^{\ast} \cL_{\alpha}] =
\chi_{c} (Z) [X, f^{\ast} \cL_{\alpha}].
$$
\end{lem}

\begin{proof}Immediate.\end{proof}

\begin{lem}Let $a$ be an integer and let
$\mu_{d} (k)$ act on $\GG_{m, k}$ by
multiplication by $\xi^{a}$, $\xi \in \mu_{d} (k)$.
For any non trivial character $\alpha$ of
$\mu_{d} (k)$, $\chi_{c} (\GG_{m, k}, \alpha) = 0$.
\end{lem}

\begin{proof}The action of $\mu_{d} (k)$ on $\GG_{m, k}$ extends to an action
on
$\PP^{1}_{k}$ leaving fixed $0$ and $\infty$. So 
it is enough to verify
that
if $\alpha$ is a non trivial character of
$\mu_{d} (k)$, $[h (\PP^{1}_{k}, \alpha)] = 0$. Now remark that, for any
$\xi$
in $k^{\times}$, the class of the graph of the multiplication by
$\xi$ in $A^{1} (\PP^{1}_{k} \times \PP^{1}_{k})$ is equal to the class
of the diagonal, hence $f_{\alpha} = 1$ if $\alpha$ is trivial,
and $f_{\alpha} = 0$ otherwise.\end{proof}

If
$f : X \rightarrow \GG_{m, k}$ and
$g : Y \rightarrow \GG_{m, k}$ are
morphisms in ${\rm Sch}_{k}$,
we denote by $f \otimes g$ the morphism
$X \times Y \rightarrow \GG_{m, k}$ given by multiplication of $f$ and $g$.
We have the following generalization of Lemma 1.4.3.

\begin{lem}Let $g : Z \rightarrow \GG_{m, k}$ be a morphism in
${\rm Sch}_{k}$.
For any character $\alpha$ of order $d$ of
$\mu_{d} (k)$ and any integer $n$ not divisible by $d$,
$$[\GG_{m, k} \times Z, ([n] \otimes g)^{\ast} \cL_{\alpha}] = 0.$$
\end{lem}

\begin{proof}Set $W = \widetilde{(\GG_{m, k} \times Z)}_{[n] \otimes g, d}$.
We may identify $W$ with
$\{(x, z, t) \in \GG_{m, k} \times Z \times \GG_{m, k} \bigm| x^{n} g(z)
= t^{d}\}$, the action of $\mu_{d} (k)$ being
multiplication on the last factor.
Set $\delta = \gcd (n, d)$, $n = \delta n'$,
$d = \delta d'$, and choose integers $a$ and $b$
such that $a n' = 1 + b d'$.
If we set $w = t^{d'} x^{- n'}$, $x' = w^{a} x$, $t' = w^{b}t$, we may identify
$W$ with
$$\{(x', z, t',w) \in \GG_{m, k} \times Z \times \GG_{m, k} \times
\GG_{m, k}
\bigm| g(z)
= w^{\delta}    \, \hbox{and} \, x'{}^{n'}
=
t'{}^{d'}
\}.$$
We may rewrite this as an isomorphism
$W \simeq \widetilde Z_{g, \delta} \times \GG_{m, k}$,
the action of $\mu_{d} (k)$ being the product of 
the action on $\widetilde Z_{g, \delta}$
given by composition with the surjection $\xi \mapsto \xi^{d'}$,
$\mu_{d} (k) \rightarrow \mu_{\delta} (k)$,
with the action on $\GG_{m, k}$ given by 
multiplication by $\xi^{a}$, for $\xi \in \mu_{d} (k)$.
By Proposition 1.3.3 (3),
$$
\chi_{c} (W, \alpha) = \sum_{\beta \in \widehat{\mu_{d} (k)}} \chi_{c}
(\widetilde Z_{g, \delta},
\beta) \,
\chi_{c} (\GG_{m, k}, \alpha \beta^{-1}).
$$
Hence, by Lemma 1.4.3 and Proposition 1.3.3 (1),
$$
\chi_{c} (W, \alpha) = \chi_{c}
(\widetilde Z_{g, \delta},
\alpha)
\chi_{c} (\GG_{m, k}).
$$
But $\chi_{c} (\widetilde Z_{g, \delta}, \alpha) = 0$, by
Proposition 1.3.3 (2), because
$\alpha$ is of order $d$ and $\delta < d$.\end{proof}

\subsection{Quotients}We discuss here motivic Euler
characteristics of quotients. This part
will only be used in
section 3.

\begin{lem}Let $X$ be a smooth projective
scheme with $G$-action, $H$ a subgroup of $G$, and
$\alpha$ a character of $G / H$. Assume
the quotient $X / H$ is smooth. Then
$h (X / H, \alpha) \simeq h (X, \alpha \circ \varrho)$,
where $\varrho$ is the projection $ G \rightarrow G / H$.
\end{lem}

\begin{proof}The projection $X \rightarrow X / H$
induces by functoriality a morphism
$$h (X / H) \longrightarrow h (X)$$
in $\cM_k$, which induces an isomorphism between
$h (X / H, \alpha)$ and $h (X, \alpha \circ \varrho)$.\end{proof}

\medskip

In view of Lemma 1.5.1 and Corollary 1.5.4, 
it seems quite natural to expect
the following statement holds.

\begin{ass}
If $X$ is 
a $G$-scheme, $H$ a subgroup of $G$, and
$\alpha$ a character of $G / H$, then
$\chi_{c} (X / H, \alpha) = \chi_{c} (X, \alpha \circ \varrho)$
and
$\chi (X / H, \alpha) = \chi (X, \alpha \circ \varrho)$,
where $\varrho$ is the projection $ G \rightarrow G / H$.
\end{ass}

In the paper \cite{Vo}, Voevodsky constructs for a perfect field $k$
a tensor triangulated category $DM_{gm} (k)$ which he calls the triangulated
category of geometrical motives. When
$k$ is of characteristic zero, he associates to any object $X$
in
${\rm Sch}_k$ 
complexes $M_{gm}^{c} (X)$ and $M_{gm} (X)$ in $DM_{gm} (k)$.
Let $E$ be a field of characteristic zero and denote by
$DM_{gm} (k)_E$ the category $DM_{gm} (k) \otimes E$.
By \cite{Vo} 2.2, when $X$ is proper and smooth
$M_{gm} (X) = M^{c}_{gm} (X)$ and the restriction of $M_{gm}$ to $\cV_k$
factorizes through an additive functor
$\cM_{k, E} \rightarrow DM_{gm} (k)_E$.
Hence there is a canonical morphism of groups
$\varphi : K_0 (\cM_{k, E}) \rightarrow K_0 (DM_{gm} (k)_E)$.
The tensor structure on $DM_{gm} (k)$ induces a ring structure on
the Grothendieck group
$K_0 (DM_{gm} (k)_E)$ and $\varphi$ is a morphism a rings.
By \cite{Vo} Corollary 3.5.5, the morphism $\varphi$
is surjective, but
it does not seem to be known whether $\varphi$
is injective or not.
It follows directly from the properties of
$\chi_c$ and $M^{c}_{gm}$ that, for any $X$ in
${\rm Sch}_k$, $\varphi (\chi_c (X)) = [M^c_{gm} (X)]$
(cf. cite{G-S} 3.2.4). Similarly, it follows by an easy
induction on dimension and the properties of
$\chi$ and $M_{gm}$ that, for any $X$ in
${\rm Sch}_k$, $\varphi (\chi (X)) = [M_{gm} (X)]$.
Let $G$ be a finite abelian group
and
assume $E$
contains all the roots of unity of order dividing $|G|$.
Let $X$ be an object of  ${\rm Sch}_{k, G}$. Then it follows
from the definition of the complexes
$M_{gm}^{c} (X)$ and $M_{gm} (X)$ that $G$
acts on them and that they decompose in $DM_{gm} (k)_E$
into direct sums of
isotypic components
$M_{gm}^{c} (X) \simeq \bigoplus_{\alpha \in \hat G} M _{gm}^{c}
(X)_{\alpha}$ and
$M_{gm} (X) \simeq \bigoplus_{\alpha \in \hat G} M_{gm} (X)_{\alpha}$.
One derives similarly as before that
$\varphi (\chi_c (X, \alpha)) = [M_{gm}^{c} (X)_{\alpha}]$
and
$\varphi (\chi (X, \alpha)) = [M_{gm} (X)_{\alpha}]$.

\begin{lem}
If $X$ is 
a $G$-scheme, $H$ a subgroup of $G$, and
$\alpha$ a character of $G / H$, then
there are canonical isomorphisms
$M_{gm}^{c} (X / H)_{\alpha} \simeq M_{gm}^{c} (X)_{\alpha \circ \varrho}$
and
$M_{gm} (X / H)_{\alpha} \simeq  M_{gm} (X)_{\alpha \circ \varrho}.$
\end{lem}

\begin{proof}This follows directly from the definition (such a statement
is already true at the level of the Nisnevich sheaves
$L^c$ and $L$ of \cite{Vo} 4.1).\end{proof}

\begin{cor}
If $X$ is 
a $G$-scheme, $H$ a subgroup of $G$, and
$\alpha$ a character of $G / H$, then
$$\varphi (\chi_{c} (X / H, \alpha))
= \varphi (\chi_{c} (X, \alpha \circ \varrho))$$
and
$$\varphi (\chi (X / H, \alpha)) = \varphi
(\chi (X, \alpha \circ \varrho)). \hfill \qed$$
\end{cor}

\section{Motivic Igusa zeta functions}
\subsection{}We will consider the ring of formal series
$K_{0} (\cM_{k, \QQ[\zeta_{d}]}) [[\LL^{-s}]]$. In this ring we will write
$\LL^{j} \cdot (\LL^{-s})^{i} = \LL^{j - s i}$, when $j \in \ZZ$ and
$i \in \NN$.
We will also consider the subring 
$K_{0} (\cM_{k, \QQ[\zeta_{d}]}) \, [\LL^{-s}]_{\rm loc}$ of
the ring
$K_{0} (\cM_{k, \QQ[\zeta_{d}]}) \, [[\LL^{-s}]]$
generated by
$K_{0} (\cM_{k, \QQ[\zeta_{d}]}) \, [\LL^{-s}]$
and the series
$$(1 - \LL^{-Ns - n})^{-1} = \sum_{i \in \NN}  \LL^{- N s i - ni},$$
for $N$ and $n$ in $\NN \setminus \{0\}$.

Let $X$ be a smooth and connected separated
$k$-scheme of finite type of
dimension $m$, $f : X \rightarrow
\AA^{1}_{k}$ be a morphism,
and $W$ be a reduced subscheme of $X$. We assume that $\mu_{d} (k)$ 
is of order $d$.
For any character $\alpha$ of $\mu_{d} (k)$ of order $d$,
we define the motivic Igusa zeta function
$\displaystyle\int_{W} (f^{s}, \alpha)$ in
$K_{0} (\cM_{k, \QQ[\zeta_{d}]}) \, [[\LL^{-s}]]$
as follows.

We denote by
$\cL_{n} (X)$ the $k$-scheme which represents the functor,
defined on the category of $k$-algebras,
$$
R \mapsto {\rm Mor}_{k-{\rm schemes}}
({\rm Spec} \, R [t]/t^{n+1} R [t], X),
$$
for $n \geq 0$ (cf. p.276 of
S. Bosch, W. L\"utkebohmert and M. Raynaud, {\it Neron models}, 
Ergeb. Math. Grenzgeb. (3) 21, Springer-Verlag, Berlin, 1990).
We denote by $\cL (X)$
the projective limit in the category of schemes of the
schemes $\cL_{n} (X)$, which exists since the transition maps are affine.
Note that for any field $K$ containing $k$ the $K$-rational points of
$\cL (X)$ are the morphisms
${\rm Spec} \, K [[t]] \rightarrow X$.
For $f : X \rightarrow
\AA^{1}_{k} = {\rm Spec} \, k [x]$ a morphism,
and $n$ in $\NN$,
we define
$Z_{n, f, W}$ as the reduced subscheme of 
$\cL (X)$
whose $K$-rational points, 
for any field $K$ containing $k$, are the morphisms $\varphi :
{\rm Spec} \, K [[t]] \rightarrow X$
sending the closed point of
${\rm Spec} \, K [[t]]$ to a point in $W$, and such that
$f \circ \varphi$ is exactly of order $n$ at the origin.
We denote by $X_{n, f, W}$ the image of $Z_{n, f, W}$
in
$\cL_{n} (X)$, viewed as a reduced subscheme of $\cL_{n} (X)$,
and by $\bar f$ the morphism
$\bar f : X_{n, f, W} \rightarrow \GG_{m, k}$
which associates
to $\varphi$ in $X_{n, f, W}$
the constant term of the series $t^{-n} x (f \circ \varphi)$.

We define,
for any character $\alpha$ of $\mu_{d} (k)$ of order $d$,
$$\int_{W} (f^{s}, \alpha)
~:=
\sum_{n \in \NN} \, 
[X_{n, f, W},  \bar f^{\ast}\cL_{\alpha}] \,\,  \LL^{-n s - (n+ 1)m}
$$
in $K_{0} (\cM_{k, \QQ[\zeta_{d}]}) \, [[\LL^{-s}]]$.
When $\alpha$ is the trivial character, we write
$\displaystyle\int_{W} f^{s}$ instead of $\displaystyle\int_{W} (f^{s},
\alpha)$.

\theoremstyle{remark}
\newtheorem*{Rks}{Remarks}
\begin{Rks}
\begin{enumerate}
\item When $\alpha$ is the trivial character,
$\displaystyle\int_{W} f^{s}$ is the image
of the series
$$\int_{W}^{\sim} f^{s} := \sum_{n \in \NN} \, 
[X_{n, f, W}] \,\,  \LL^{-n s - (n+ 1)m}
$$ in 
$K_{0} ({\rm Sch}_{k}) [\LL^{-1}]  [[\LL^{-s}]]$ by the natural
morphism
$$K_{0} ({\rm Sch}_{k}) [\LL^{-1}]  [[\LL^{-s}]] \longrightarrow
K_{0} (\cM_{k}) \, [[\LL^{-s}]]$$ induced by $\chi_{c}$.
\item It would be interesting to investigate whether motivic Igusa
functions already exist at a finer level than a Grothendieck group
of Chow motives, for instance at the level
of
complexes of Chow motives, or objects of $DM_{gm} (k)$, or ``mixed motives''.
\end{enumerate}
\end{Rks}
\medskip

\subsection{}Let $D$ be the divisor defined by $f = 0$ in
$X$. Let $(Y, h)$ be a resolution of $f$. By this, we mean that 
$Y$ is a smooth and connected $k$-scheme of finite type, $h : Y \rightarrow
X$
is proper, that the restriction
$h : Y \setminus h^{-1} (D) \rightarrow
X \setminus D$ is an isomorphism, and that
$(h^{-1} (D))_{\rm red}$ has only normal crossings as a
subscheme
of $Y$. Let $E_{i}$, $i \in J$, be the irreducible (smooth) components
of
$(h^{-1} (D))_{\rm red}$. For each $i \in J$, denote by
$N_{i}$ the multiplicity of $E_{i}$ in the divisor of
$f \circ h$ on $Y$, and by
$\nu_{i} - 1$ the multiplicity of $E_{i}$ in the divisor
of $h^{\ast} dx$, where $d x$ is a local non vanishing volume form,
{\it i.e.} a local generator of the sheaf of differential forms of
maximal degree.
For $i \in J$ and $I \subset J$, we consider the schemes
$E_{i}^{\circ} := E_{i} \setminus \cup_{j \not= i} E_{j}$,
$E_{I} := \cap_{i \in I} E_{i}$, and
$E_{I}^{\circ} := E_{I} \setminus \cup_{j \in J \setminus I} E_{j}$.
When $I =\emptyset$, we have $E_{\emptyset} = Y$.

Now denote by $J_{d}$ the set of $I \subset J$ such that
$d \, \vert \, N_{i}$ for all $i $ in $I$
and by
$U_{d}$ the union of the $E_{I}^{\circ}$,
with $I$ in $J_{d}$.
Let $Z$ be locally closed in $U_{d}$.
For any character $\alpha$ of $\mu_{d} (k)$ of order $d$,
we will construct an element $[Z_{f, \alpha}]$ in
$K_{0} (\cM_{k, \QQ[\zeta_{d}]})$
as follows. If on $Z$ we
may write $f \circ h = u v^{d}$ with $u$ non vanishing on 
$Z$, we set
$[Z_{f, \alpha}] = [Z, u^{\ast} \cL_{\alpha}]$.  It  is well defined by
Lemma 1.4.1. In general we cover $Z$
by a finite set of $Z_{r}$'s for which the previous condition holds,
and we set
$$
[Z_{f, \alpha}] = \sum_{r} [(Z_{r})_{f, \alpha}]
-
\sum_{r_{1} \not = r_{2}} [(Z_{r_{1}} \cap  Z_{r_{2}})_{f, \alpha}]
+ \cdots,
$$
which is well defined by additivity of $\chi_{c} (\quad, \alpha)$.

\medskip

We can now state the following result.

\begin{theorem}For any character $\alpha$ of $\mu_{d} (k)$ of order $d$,
$$\int_{W} (f^{s}, \alpha)
=
\LL^{-m}
\sum_{I \in J_{d}}
[(E_{I}^{\circ} \cap h^{-1} (W))_{f, \alpha}] \prod_{i \in I}
\frac
{(\LL - 1) \, \LL^{-N_{i}s - \nu_{i}}}
{1 - \LL^{-N_{i}s - \nu_{i}}}
$$
in $K_{0} (\cM_{k, \QQ[\zeta_{d}]}) \, [[\LL^{-s}]]$.
In particular $\displaystyle\int_{W} (f^{s}, \alpha)$ belongs to
the ring
$K_{0} (\cM_{k, \QQ[\zeta_{d}]}) \, [\LL^{-s}]_{\rm loc}$.
\end{theorem}

\begin{proof}We set $\cL (X, D) := \cL (X) \setminus \cL (D)$, and we
define
similarly the scheme
$\cL (Y, h^{-1} (D))$.
We denote by $\pi_{n}$ the projections
$\cL (X, D) \rightarrow \cL_{n} (X)$ and
$\cL (Y, h^{-1} (D)) \rightarrow \cL_{n} (Y)$, and by $\pi$ the
projections onto
$X$ and $Y$ respectively.
If $U$ is a reduced subscheme of $Y$ we set
$\cL_{U} (Y, h^{-1} (D)) = \pi^{-1} (U)$
and
$\cL_{U} (X, D) = h_{\ast} (\cL_{U} (Y, h^{-1} (D)))$,
and we define similarly 
$\cL_{n, U} (Y)$.
Moreover, for $n' \geq n$, we denote by $\pi_{n', n}$ the projections
$\cL_{n'} (X) \rightarrow \cL_{n} (X)$ and 
$\cL_{n'} (Y) \rightarrow \cL_{n} (Y)$.

The morphism $h$ being proper, composition with $h$ induces a
bijective morphism
$h_{\ast} : \cL (Y, h^{-1} (D)) \simeq \cL (X, D)$,
and we have a commutative diagram 
\begin{equation*}\xymatrix{
\cL (Y, h^{-1} (D)) \ar[r]^<<<<<{\pi_{n}} \ar[d]^{h_{\ast}}& \cL_{n} (Y)
\ar[d]^{h_{n \ast}}\\
\cL (X, D) \ar[r]^{\pi_{n}}& \cL_{n} (X).
}
\end{equation*}

For $U$ a reduced subscheme of $Y$,
we set 
$$X_{n, f, W, U}^{n'} := \pi_{n'} (Z_{n, f, W} \cap  \cL_{U} (X, D)).$$
We define
$$
M_{n, f, W, U} :=
\LL^{- (n - n')m} \,
[X_{n, f, W, U}^{n'}, (\bar f \circ \pi_{n', n})^{\ast} \cL_{\alpha}],
$$
where $n'$ is big enough with respect to $n$.
By Lemma 1.4.2, the definition of
$M_{n, f, W, U}$ does not depend on $n'$ because $\pi_{n', n}$ is a
locally trivial fibration for the Zariski topology with fiber
$\AA_{k}^{m (n' - n)}$ and because
$Z_{n, f, W} \cap  \cL_{U} (X, D)$
is a union of fibers of
$\pi_{n'} : \cL (X, D) \rightarrow \cL_{n'} (X)$
when $n' \gg n$, since $h_{\ast}^{-1}$ has ``only powers of $f$ in the
denominator''.

For any character $\alpha$ of $\mu_{d} (k)$ of order $d$, we define
$$\int_{U \cap h^{-1} (W)}  h^{\ast}(f^{s}, \alpha)
=
\sum_{n \in \NN} \, 
M_{n, f, W, U} \,\,  \LL^{-n s - (n+ 1)m}
$$
in $K_{0} (\cM_{k, \QQ[\zeta_{d}]}) \, [[\LL^{-s}]]$.
The result is a direct consequence of the following proposition,
by additivity of $\chi_{c} (\quad, \alpha)$.\end{proof}

\begin{prop}Assume $U \subset E^{\circ}_{I}$ and
$f \circ h = u \prod_{i \in I} y_{i}^{N_{i}}$ on a neighborhood of
$U$,
where $u$ is a unit on $U$, and $y_{i} = 0$ is an equation for
$E_{i}$ on a neighbourhood of
$U$.
\begin{enumerate}
\item[(1)]If $d$ divides $N_{i}$, for all $i \in I$,
then
$$\int_{U \cap h^{-1} (W)}  h^{\ast}(f^{s}, \alpha)
=
\LL^{-m}
[(U\cap h^{-1} (W))_{f, \alpha}] \prod_{i \in I}
\frac
{(\LL - 1) \, \LL^{-N_{i}s - \nu_{i}}}
{1 - \LL^{-N_{i}s - \nu_{i}}}
$$
in $K_{0} (\cM_{k, \QQ[\zeta_{d}]}) \, [[\LL^{-s}]]$.
\item[(2)]If $d$ does not divide $N_{i}$, for some $i \in I$,
then
$$\int_{U \cap h^{-1} (W)}  h^{\ast}(f^{s}, \alpha) = 0.$$
\end{enumerate}
\end{prop}

\begin{proof}We may from the beginning assume $U \cap h^{-1} (W) = U$
and we will write $X_{n, f, U}^{n'}$ instead of
$X_{n, f, W, U}^{n'}$. We will use the following lemma.

\bigskip

Let $X$, $Y$ and $F$ be algebraic varieties over $k$,
and let
$A$, {\it resp.} $B$, be a constructible subset of $X$,
{\it resp.} $Y$. We say that
a map
$\pi : A \rightarrow B$ is {\it piecewise trivial fibration with fiber}
$F$, if there exists a finite partition of $B$ in subsets $S$ which are
locally closed
in $Y$ such that $\pi^{- 1} (S)$ is locally closed in $X$ and
isomorphic, as
a variety over $k$, to $S \times F$, with $\pi$ corresponding
under the isomorphism to the projection
$S \times F \rightarrow S$. We say that the map $\pi$ is
a
{\it piecewise trivial fibration over} some constructible subset $C$ of
$B$,
if the restriction of $\pi$ to $\pi^{- 1} (C)$ is a piecewise 
trivial fibration.

\begin{lem}Let $X$ and $Y$ be connected smooth schemes over
a field $k$ and let $h : Y \rightarrow X$ be a birational morphism.
For $e$ in $\NN$, let $\Delta_e$ be the reduced
subscheme of $\cL (Y)$ defined
by
$$
\Delta_e (K)
:= \{\varphi \in Y (K[[t]]) \bigm \vert {\rm ord}_t {\rm det} \cJ_{\varphi}
= e \},
$$
for any field $K$ containing $k$,
where $\cJ_{\varphi}$ is the jacobian of $h$ at $\varphi$.
For $n$ in $\NN$, let
$h_{n \ast} : \cL_n (Y) \rightarrow \cL_n (X)$ be the morphism
induced by $h$, and let $\Delta_{e, n}$ be the image of
$\Delta_e$ in $\cL_n (Y)$. If $n \geq 2e$, the following holds.
\begin{enumerate}
\item[a)]The set $\Delta_{e, n}$ is  a union of fibers of
$h_{n \ast}$.
\item[b)]The restriction of $h_{n \ast}$ to $\Delta_{e, n}$
is a piecewise trivial fibration
with fiber
$\AA^e_k$ onto its image.
\end{enumerate}
\end{lem}

\begin{proof}This is a special case of Lemma 3.4 of \cite{D-L2}.\end{proof}

\medskip

Let $m_{i}$, $i \in I$, be strictly positive
integers with $\sum_{i \in I} m_{i} N_{i} = n$.
We denote by $\widetilde Z_{(m_{i}), U}$
the reduced subscheme of $\cL_{U} (Y, h^{-1} (D))$
whose $K$-rational points $\varphi :
{\rm Spec} \, K [[t]] \rightarrow Y$,
for any field $K$ containing $k$, satisfy the condition that
$y_{i} \circ \varphi$ is exactly of order $m_{i}$ at the origin,
for $i \in I$. We
denote by $Y_{(m_{i}), U}$ the image of $\widetilde Z_{(m_{i}), U}$
in
$\cL_{n, U} (Y)$.
By Lemma 2.2.3, for $n'$ big enough with respect to $n$, the set
$X^{n'}_{n, f, U}$ is the disjoint finite
union of the sets
$h_{n' \ast}(\cup_{\cE_{e}} \, \pi_{n', n}^{-1} (Y_{(m_{i}), U}))$
for $e = 0, 1, 2, \ldots$,
where $\cE_{e}$ is the set of all $(m_{i})_{i \in I}$ with
$m_{i} > 0$, $\sum_{i \in I} m_{i} N_{i} = n$
and $\sum_{i \in I} (\nu_{i} - 1) m_{i} = e$.
Hence we deduce from Lemma 1.4.2 and Lemma 2.2.3
\begin{multline*}
\int_{U \cap h^{-1} (W)}  h^{\ast}(f^{s}, \alpha)
=\\
\sum_{m_{i} > 0} \LL^{- \sum_{i \in I} (\nu_{i} - 1) m_{i}}
[Y_{(m_{i}), U}, (\bar f \circ h_{n \ast})^{\ast} \cL_{\alpha}] \,
\LL^{- ns - (n + 1) m},
\end{multline*}
with $n = \sum_{i \in I} m_i N_i$. (Actually we need here the slightly
stronger
version of Lemma 2.2.3 obtained by replacing $\Delta_{e}$ by
$\Delta_{e} \cap \cL_{U} (Y, h^{-1} (D))$. But the proof of this version
is the same.)

Now remark that $\gamma : Y_{(m_{i}), U} \rightarrow U$ is a
locally trivial fibration for the Zari\-ski topology
with fibre $\GG_{m}^{I} \times \AA^{n m - \sum_{i \in I} m_i}$.
On $Y_{(m_{i}), U}$ the function $\bar f \circ h_{n \ast}$
coincides with the product 
$(u_{\vert U} \circ \gamma) \cdot \psi$, with $\psi (\varphi)$
is the constant term of $[\prod_{i \in I} y_{i} (\varphi (t))^{N_{i}}]
t^{- \sum_{i} m_{i} N_{i}}$.

So, if $d$ divides $N_{i}$ for all $i \in I$,
we deduce from Lemma 1.4.1 and Lemma 1.4.2
that
$$[Y_{(m_{i}), U}, (\bar f \circ h_{n \ast})^{\ast} \cL_{\alpha}]
=
[U_{f, \alpha}] (\LL - 1)^{\vert I \vert} \, \LL^{n m - \sum_{i \in I} m_i},
$$
and the result follows from the previous relation.

Assume now that some $N_{i}$ with $i \in I$
is not divisible by
$d$. We may assume, shrinking $U$
if necessary, that $\gamma$ is a product. We may then identify
$Y_{(m_{i}), U}$ with a product $\GG_{m, k} \times Z$ in such a way
that $\bar f \circ h_{n \ast} = [m] \otimes g$, with $m$ 
not divisible by $d$
(notations of 1.4)
and now the result follows from Lemma 1.4.4.\end{proof}

\medskip

\subsection{Relation with the topological zeta functions of \protect\cite{D-L1}}Let us denote
by
$K_{0} (\cM_{k, \QQ[\zeta_{d}]}) \, [\LL^{-s}]'_{\rm loc}$ the subring
of
$K_{0} (\cM_{k, \QQ[\zeta_{d}]}) \, [\LL^{-s}]_{\rm loc}$
generated by the ring of polynomials
$K_{0} (\cM_{k, \QQ[\zeta_{d}]}) \, [\LL^{-s}]$
and by the quotients
$(\LL - 1) (1 - \LL^{-Ns - n})^{-1}$,
for $N$ and $n$ in $\NN \setminus \{0\}$.
By expanding $\LL^{-s}$ and $(\LL - 1) (1 - \LL^{-Ns - n})^{-1}$ into series
in $\LL - 1$, one gets a canonical morphism of algebras
\begin{multline*}
\varphi : K_{0} (\cM_{k, \QQ[\zeta_{d}]}) \,  [\LL^{-s}]'_{\rm loc}
\longrightarrow\\
{\overline K}_{0} (\cM_{k, \QQ [\zeta_{d}]}) \, [s]
[(N s + n)^{-1}]_{n, N \in \NN \setminus \{0\}}
[[\LL - 1]],
\end{multline*}
where $[[\LL - 1]]$ denotes completion with respect to the ideal
generated by $\LL - 1$ and where
${\overline K}_{0} (\cM_{k, \QQ [\zeta_{d}]})$ is the largest quotient of
$K_{0} (\cM_{k, \QQ [\zeta_{d}]})$ with no $(\LL - 1)$-torsion,
cf. remark 1.2.3.
Taking the
quotient of
${\overline K}_{0} (\cM_{k, \QQ [\zeta_{d}]}) \, [s]
[(N s + n)^{-1}]_{n, N \in \NN \setminus \{0\}}
[[\LL - 1]]$
by the ideal generated by
$\LL - 1$,
one obtains the evaluation morphism 
\begin{multline*}
{\rm ev}_{\LL = 1} :
{\overline K}_{0} (\cM_{k, \QQ [\zeta_{d}]}) \, [s]
[(N s + n)^{-1}]_{n, N \in \NN \setminus \{0\}}
[[\LL - 1]]
\longrightarrow \\
 ({\overline K}_{0} (\cM_{k, \QQ [\zeta_{d}]}) / \LL - 1) \, [s]
[(N s + n)^{-1}]_{n, N \in \NN \setminus \{0\}}.
\end{multline*}

For $X$ in $\cV_{k}$, we denote by $\chi_{\rm top} (X)$ the usual
Euler characteristic of $X$
(say in \'{e}tale $\bar \QQ_{\ell}$-cohomology).
This induces by 1.2.3 a morphism
$$\chi_{\rm top} : 
{\overline K}_{0} (\cM_{k, \QQ[\zeta_{d}]})
\longrightarrow \ZZ,$$ which induces, since
$\chi_{\rm top} (\LL) = 1$, 
a morphism  
\begin{multline*}\chi_{\rm top} :
({\overline K}_{0} (\cM_{k, \QQ [\zeta_{d}]})  / \LL - 1) \, [s] 
[(N s + n)^{-1}]_{n, N \in \NN \setminus \{0\}}\\
\longrightarrow  \ZZ
[s]
[(N s + n)^{-1}]_{n, N \in \NN \setminus \{0\}}.
\end{multline*}

By Theorem 2.2.1,
for any character $\alpha$ of $\mu_{d} (k)$ of order $d$, the motivic Igusa function
$\displaystyle\int_{W} (f^{s},
\alpha)$ belongs to $K_{0} (\cM_{k, \QQ[\zeta_{d}]}) \, [\LL^{-s}]'_{\rm loc}$,
hence we can consider
the rational function $(\chi_{\rm top} \circ {\rm ev}_{\LL = 1} \circ \varphi)
(\displaystyle\int_{W} (f^{s}, \alpha))$ in $\CC (s)$.

\begin{prop}For any character $\alpha$ of $\mu_{d} (k)$ of order $d$,
$$(\chi_{\rm top} \circ {\rm ev}_{\LL = 1} \circ \varphi)
\Bigl(\int_{W} (f^{s}, \alpha)\Bigr)
=
\sum_{I \in J_{d}}
\chi_{\rm top} (E_{I}^{\circ} \cap h^{-1} (W)) \prod_{i \in I} \,
\frac
{1}{N_{i}s + \nu_{i}}.
$$
\end{prop}

\begin{proof}By Theorem 2.2.1, it is enough to check that
$$\chi_{\rm top} ([(E_{I}^{\circ} \cap h^{-1} (W))_{f, \alpha}])
=
\chi_{\rm top} (E_{I}^{\circ} \cap h^{-1} (W)).$$ This is clear, since the
\'{e}tale $\ell$-adic realization of $[(E_{I}^{\circ} \cap h^{-1} (W))_{f, \alpha}]$
is given by cohomology with compact support of a rank one lisse sheaf on
$E_{I}^{\circ} \cap h^{-1} (W)$.\end{proof}

\begin{Rks}
\begin{enumerate}
\item It follows from Proposition 2.3.1 that the topological zeta
functions
of \cite{D-L1} are obtained by specialization of motivic Igusa zeta functions. This gives another
proof, not using $p$-adic analysis, of the main results 
of  \cite{D-L1} on the invariance of topological zeta functions
in the algebraic case. In fact
it is easily checked that similar arguments work also in the complex
analytic case.
\item It might be interesting to study the functions
$({\rm ev}_{\LL = 1} \circ \varphi)
(\displaystyle\int_{W} (f^{s}, \alpha))$ which belong to
$({\overline K}_{0} (\cM_{k, \QQ [\zeta_{d}]}) / \LL - 1) \, (s)$.
\end{enumerate}
\end{Rks}

\subsection{Relation with $p$-adic Igusa local zeta functions}Let $p$
be a prime number and let $K$ be a finite extension of $\QQ_p$.
Let $R$ be the valuation ring of $K$, $P$ the maximal ideal of $R$, and
$\bar K = R / P$ the residue field of $K$. Let $q$ denote the cardinality of
$\bar K$, so $\bar K \simeq \FF_q$. For $z$ in $K$,  let
${\rm ord} \, z$ denote the valuation of $z$, and  set
$\vert z \vert = q^{- {\rm ord} \, z}$ and ${\rm ac} (z) = z \pi^{- {\rm ord} \, z}$,
where $\pi$ is a fixed uniformizing parameter of $R$.
Let $f$ be an element of $R [x_1, \ldots, x_m]$ which is not zero
modulo $P$.
For any character $\alpha : R^{\times} \rightarrow \CC^{\times}$,
one defines the $p$-adic Igusa local zeta function $Z (s, \alpha)$
as the integral
$$
Z (s, \alpha) = \int_{R^m} \alpha ({\rm ac} (f (x))) \vert f (x) \vert^s \vert dx \vert,
$$
for $s \in \CC$, ${\rm Re} (s) > 0$, 
where $\vert dx \vert$ denotes the Haar measure on
$K^m$ normalized in such of way that $R^m$ is of volume 1.

Let $(Y, h)$ be a resolution of $f$ as in 2.2. We say the resolution $(Y, h)$
has good reduction ${\rm mod} \, P$, if $Y$ has a smooth 
model $Y_R$
over ${\rm Spec} \, R$ such that $h$ extends to a morphism
$Y_{R} \rightarrow \AA^{m}_{R}$ and 
such that the closure of $h^{- 1} (D)_{\rm red}$ in $Y_R$
is a relative divisor with normal crossings
over $Y_R$. For $Z$ closed in
$Y$, we denote by $\overline Z$ the fiber over the closed point
of the closure of $Z$ in
$Y_R$. Hence
${\overline Y}$ and all the ${\overline E}_i$'s are smooth,
$\cup_{i \in J} {\overline E}_i$ is a divisor with normal crossings, and the schemes
${\overline E}_i$ and
${\overline E}_j$ have no component in common for $i \not= j$. Let
$(Y, h)$ be a resolution with good reduction ${\rm mod} \, P$. For $I \subset J$
we
have $\overline{E}_I = \cap_{i \in I} {\overline E}_i$ and we set
$\overline{E}_{I}^{\circ} := \overline{E}_{I}
\setminus \cup_{j \in J \setminus I} \overline{E}_{j}$.

Assume now the character $\alpha$ is of finite order $d$
and is trivial on $1 + P$.
Choose a prime number $\ell \not= p$ and denote by $\cL_{\alpha}$ the
Kummer $\bar \QQ_{\ell}$-sheaf on $\GG_{m, \bar K}$ associated to $\alpha$
viewed as a character of ${\bar K}^{\times}$ (here we choose an
isomorphism between the group of roots of unity in $\CC$
and in $\bar \QQ_{\ell}$).
Set $U = {\overline Y} \setminus \cup_{i \in J} {\overline E}_i$
and denote by $\nu : U \hookrightarrow {\overline Y}$ the open immersion
and by $\beta : U \rightarrow \GG_{m, \bar K}$ the map induced
by ${\overline f} \circ {\overline h}$.
Set $\cF_{\alpha} = \nu_{\ast} \beta^{\ast} \cL_{\alpha}$.
Denote by ${\overline K}^a$ the algebraic closure
of ${\overline K}$  and by $F$ the geometric Frobenius automorphism.

In the good reduction case the following result gives a cohomological
expression for $p$-adic Igusa local zeta functions.

\begin{theorem}[\protect\cite{D1}\protect\cite{D2}]Let
$(Y, h)$ be a resolution of $f$
with good reduction ${\rm mod} \, P$. Assume the character
$\alpha$ is of finite order $d$
and is trivial on $1 + P$.
Then
$$
Z (s, \alpha) = q^{- m}
\sum_{I \in J_{d}}
c_{I, \alpha}
\prod_{i \in I}
\frac
{(q - 1) \, q^{-N_{i}s - \nu_{i}}}
{1 - q^{-N_{i}s - \nu_{i}}}
$$
with
$$
c_{I, \alpha} = \sum_i (- 1)^i {\rm Tr} \, (F,
H^i_c ({\overline E}^{\circ}_I \otimes
{\overline K}^a, \cF_{\alpha})).
$$
\end{theorem}

In conclusion, in view of Theorem 2.2.1 and Theorem 2.4.1,
one can state that
``in the good reduction case, the $p$-adic Igusa local zeta functions
are given by the trace of the Frobenius action on the $\ell$-adic \'{e}tale realization of the
corresponding motivic ones''.

As in the $p$-adic case (see e.g. \cite{D3}), there is the intriguing
question whether $\int_{X} (f^{s}, \alpha)$ always belong to
$K_{0} (\cM_{k, \QQ [\zeta_{d}]})
[(1 - \LL^{-Ns - n})^{-1}]_{(N, n) \in M},$
where $M$ is the set of all pairs $(N, n)$ in $(\NN \setminus
\{0\})^{2}$
with $\exp (2 \pi i n / N)$ an eigenvalue of the monodromy action on the
complex $R \psi_{f}$ of
nearby cycles on $f^{-1} (0)$. For some recent
work in the $p$-adic case, see \cite{Veys1}, \cite{Veys2}.


\subsection{Relation with motivic integration}
M.Kontsevich introduced in \cite{K} the completion
$\widehat K_{0} ({\rm Sch}_{k})$ of 
$K_{0} ({\rm Sch}_{k}) [\LL^{-1}]$ with respect to the filtration
$F^m K_{0} ({\rm Sch}_{k}) [\LL^{-1}]$, where $F^m
K_{0} ({\rm Sch}_{k}) [\LL^{-1}]$ is
the subgroup of $K_{0} ({\rm Sch}_{k}) [\LL^{-1}]$ generated by
$\{ [S] \, \LL^{- i} \bigm \vert
i - \dim S \geq m\}$,
and defined, for smooth $X$ over $k$, a motivic integration on $\cL (X)$ with
values
into ${\widehat K_{0} ({\rm Sch}_{k})}$. In the paper \cite{D-L2}, we
extended Kontsevich's construction to semi-algebraic subsets of $\cL
(X)$ and also to the non smooth
case. 
The following statement is proved in \cite{D-L2} (Definition-Proposition 3.2).

\begin{df}Let $X$ be an algebraic
variety
over $k$ of pure dimension $m$. Denote by
$\pi_{n}$ the natural morphism
$\cL (X) \rightarrow \cL_{n} (X)$.
Let $\BB$ be the boolean algebra of all semi-algebraic
subsets
of $\cL (X)$. There exists a unique map $\mu : \BB \rightarrow
\widehat K_{0} ({\rm Sch}_{k})$ satisfying the following three properties.
\begin{enumerate}
\item[] (2.5.2) \, If $A \in
\BB$ is stable at level $n$, then
$$\mu (A) = [\pi_{n} (A)] \LL^{- (n + 1) m}.$$
\item[] (2.5.3) \, If $A \in
\BB$ is contained in $\cL (S)$ with $S$ a closed subvariety of $X$ with 
${\rm dim} \, S < {\rm dim} \, X$, then $\mu (A) = 0$.
\item[] (2.5.4) \, Let $A_{i}$ be in $\BB$ for each $i$ in $\NN$.
Assume that the
$A_{i}$'s are mutually disjoint and that
$A := \bigcup_{i \in \NN} A_{i}$ is semi-algebraic. Then
$\sum_{i \in \NN} \mu (A_{i})$ converges in $\widehat K_{0} ({\rm Sch}_{k})$
to $\mu (A)$.
\end{enumerate}
We call this unique map the {\it motivic volume} on $\cL (X)$
and denote it by $\mu_{\cL (X)}$ or $\mu$. Moreover we have
\begin{enumerate}
\item[] (2.5.5) \, If $A$ and $B$ are in $\BB$, $A \subset B$,
and if $\mu (B)$ belongs to the closure $F^{m} (\widehat K_{0} ({\rm Sch}_{k}))$
of $F^{m}K_{0} ({\rm Sch}_{k}) [\LL^{-1}]$ in $\widehat K_{0} ({\rm Sch}_{k})$,
then $\mu (A) \in F^{m} (\widehat K_{0} ({\rm Sch}_{k}))$.
\end{enumerate}
Hence, for $A$ in $\BB$ and $\alpha : A \rightarrow \ZZ \cup
\{+ \infty\}$ a simple
function,
we can define
$$
\int_{A} \LL^{- \alpha} d \mu := \sum_{n \in \ZZ} \mu (A \cap \alpha^{-1} (n))
\, \LL^{- n}
$$
in $\widehat K_{0} ({\rm Sch}_{k})$, whenever the right hand side
converges in $\widehat K_{0} ({\rm Sch}_{k})$, in which case we say that $\LL^{- \alpha}$
is integrable on $A$. If the function
$\alpha$ is bounded from below, then
$\LL^{- \alpha}$
is integrable on $A$, because of (2.5.5).
\end{df}

Semi-algebraic
subsets
of $\cL (X)$ and simple functions on semi-alge\-bra\-ic
subsets are defined in \cite{D-L2}, as well as the notion of
stable semi-algebraic
subsets
of $\cL (X)$ of level $n$. In particular,
$\cL (X)$ is a semi-algebraic subset of $\cL (X)$
and, for any morphism $g : Y \rightarrow X$ of algebraic varieties over
$k$,
the image of $\cL (Y)$ in $\cL (X)$ under the morphism induced by $g$ is
a semi-algebraic subset of $\cL (X)$.
When $X$ is smooth, a 
semi-algebraic subset of $\cL (X)$
is stable of level $n$ if and only if it is a union of fibers
of $\pi_{n} : \cL (X) \rightarrow \cL_{n} (X)$.
Consider a coherent sheaf of ideals $\cI$
on $X$
and  denote by 
${\rm ord}_t \cI$
the function
${\rm ord}_t \cI : \cL (X) \rightarrow \NN \cup \{+\infty\}$ given by
$\varphi \mapsto \min_{g} {\rm ord}_t g ( \varphi)$,
where the minimum is taken over all $g$ in the stalk $\cI_{\pi_{0}
(\varphi)}$
of $\cI$ at $\pi_{0}
(\varphi)$. The function ${\rm ord}_t \cI$ is a simple function.
When $X$ is smooth and $\cI$ is the ideal sheaf of an effective
divisor $D$ on $X$, the motivic integral
$\int_{\cL (X)} \LL^{- {\rm ord}_{t} \cI} d \mu$ was first introduced by
Kontsevich \cite{K} and denoted by him $[\int_{X} e^{D}]$.
In particular, for a morphism $f : X \rightarrow \AA^{1}_{k}$ with
divisor
$D$ and
a natural number $d$ in $\NN$,
we
can consider the
motivic integral
$\int_{\pi_{0}^{-1}(W)} \LL^{- {\rm ord}_{t} \cO (-d D)} d \mu$,
for any reduced subscheme $W$ of $X$, because $\pi_{0}^{-1} (W)$ is a
semi-algebraic subset of $\cL (X)$.
It follows from Theorem 5.1 of \cite{D-L2} that
$\int_{\pi_{0}^{-1}(W)} \LL^{- {\rm ord}_{t} \cO (-d D)} d \mu$
belongs to the image of
$K_{0} ({\rm Sch}_{k}) [\LL^{-1}, ((\LL^{i} - 1)^{-1})_{i
\geq 1}]$ in $\widehat K_{0} ({\rm Sch}_{k})$.
On the other hand, the motivic Igusa function
$\int_{W} f^{s}$ belongs to
$K_{0} (\cM_{k}) \, [\LL^{-s}]_{\rm loc}$
and is the natural image of a well defined element
$\int_{W}^{\sim} f^{s}$ in $K_{0} ({\rm Sch}_{k})
[\LL^{-1}] [[\LL^{-s}]]$, cf. remark 1 in 2.1. Moreover the
proof of Theorem 2.2.1 also shows that
$\int_{W}^{\sim} f^{s}$ belongs to  $K_{0} ({\rm Sch}_{k})
[\LL^{-1}] [\LL^{-s}]_{\rm loc}$.
Hence
for any 
natural number $d$ in $\NN$, we can formally replace $s$ by $d$
and obtain by evaluation
a well defined element $(\int_{W}^{\sim} f^{s})_{\vert s = d}$
in $\widehat K_{0} ({\rm Sch}_{k})$.
The following statement is a direct
consequence of the definitions.

\setcounter{theorem}{5}
\begin{prop}Let $X$ be a smooth and connected separated
$k$-scheme of
finite type of pure
dimension $m$, $f : X \rightarrow
\AA^{1}_{k}$ be a morphism,
and $W$ be a reduced subscheme of $X$. For any 
natural number $d$ in $\NN$, the equality
$$
\int_{\pi_{0}^{-1}(W)} \LL^{- {\rm ord}_{t} \cO (-d D)} d \mu =
\Bigl(\int_{W}^{\sim} f^{s}\Bigr)_{\vert s = d}
$$
holds in $\widehat K_{0} ({\rm Sch}_{k})$.\hfill$\qed$
\end{prop}

\section{Functional equation}
\subsection{}We denote by 
$K_{0} (\cM_{k, \QQ[\zeta_{d}]}) \, [\LL^{s}, \LL^{- s}]_{\rm loc}$
the localization of the algebra of Laurent polynomials
$K_{0} (\cM_{k, \QQ[\zeta_{d}]}) \, [\LL^{s}, \LL^{- s}]$
with respect to the multiplicative set generated by the polynomials
$1 - \LL^{-Ns - n}$,
for $N$ and $n$ in $\NN \setminus
\{0\}$.
One may consider
$K_{0} (\cM_{k, \QQ[\zeta_{d}]}) \, [\LL^{- s}]_{\rm loc}$ as embedded in
$K_{0} (\cM_{k, \QQ[\zeta_{d}]}) \, [\LL^{s}, \LL^{- s}]_{\rm loc}$.
The involution
$M \mapsto M^{\vee}$ extends
to $K_{0} (\cM_{k, \QQ[\zeta_{d}]})$. One can extend it
to a $K_{0} (\cM_{k, \QQ[\zeta_{d}]})$-algebra involution
on
$K_{0} (\cM_{k, \QQ[\zeta_{d}]}) \, [\LL^{s}, \LL^{-s}]_{\rm loc}$
by setting $(\LL^{s})^{\vee} = \LL^{-s}$,
$(\LL^{ -s})^{\vee} = \LL^{s}$,
and $((1 - \LL^{-Ns - n})^{-1})^{\vee} = - \LL^{-Ns - n}
(1 - \LL^{-Ns - n})^{-1}$.

\medskip
In this section we assume $X = \AA^{m}_{k}$
and $f$ is a homogenous polynomial
of degree $r$.

\begin{theorem}\begin{enumerate}
\item[(1)]The equality
$$
\Bigl(\int_{X} f^{s} \Bigr)^{\vee}
=
\LL^{- rs}
\int_{X} f^{s}
$$
holds in
$K_{0} (\cM_{k}) \, [\LL^{s}, \LL^{-s}]_{\rm loc}$.
\item[(2)]Assume 1.5.2 holds.
Then, for any character $\alpha$ of $\mu_{d} (k)$ of order $d$,
$$\Bigl(\int_{X} (f^{s}, \alpha)\Bigr)^{\vee}
=
\LL^{- rs}
\int_{X} (f^{s}, \alpha^{-1})
$$
in
$K_{0} (\cM_{k, \QQ[\zeta_{d}]}) \, [\LL^{s}, \LL^{-s}]_{\rm loc}$.
\end{enumerate}
\end{theorem}

We  begin with the following lemma.

\begin{lem}Let $\alpha$ be a character
of $\mu_{d} (k)$ of order $d$. If $d$ does not divide
$r$, then
$\displaystyle\int_{X} (f^{s}, \alpha) = 0.$
\end{lem}

\begin{proof}
It is enough to prove that
$[X_{n, f, \AA^{m}_{k}}, \bar f^{\ast} \cL_{\alpha}] = 0$ if
$d$ does not divide
$r$. From the second displayed formula in the proof of Proposition 3.2.1 below,
which actually holds for any $d$, it follows that it suffices to prove
that $[X_{n, f, \{ 0\}}, \bar f^{\ast} \cL_{\alpha}] = 0$. Thus we
have to show that
$\displaystyle\int_{\{ 0\}} (f^{s}, \alpha) = 0.$
Let $\bar D$, {\it resp.} $D$, be the divisor in $\PP^{m - 1}_{k}$,
{\it resp.} $\AA^{m}_{k}$, defined by $f = 0$, and let $\bar h :
\bar Y \rightarrow \PP^{m - 1}_{k}$ be a resolution of
$\bar D \subset \PP^{m - 1}_{k}$ (in the sense of 2.2).
Denote by $\pi : B \rightarrow \AA^{m}_{k}$
the blowing up of $\{0\}$ in $\AA^{m}_{k}$, and by $p$
the natural map $p : B \rightarrow \PP^{m - 1}_{k}$
which is the identity on $\pi^{-1} \{0\} = \PP^{m - 1}_{k}$. Note
that $p$ is a locally trivial fibration for the Zariski topology with
fiber $\AA^{1}_{k}$. One verifies that the natural map
$$
h : \bar Y \times_{\PP^{m - 1}_{k}} B \rightarrow
B \rightarrow \AA^{m}_{k} 
$$
is a resolution of $D \subset \AA^{m}_{k}$. Moreover $h^{-1} (D)
\simeq \bar Y$ is a component (hence equal to some $E_{i}$) of $h^{-1}
(D)$ on which $f \circ h$ has multiplicity $r$. Thus 
$\displaystyle\int_{\{ 0\}} (f^{s}, \alpha) = 0$ when $d$ does not
divide $r$, by Theorem 2.2.1.
\end{proof}

\subsection{Proof of Theorem 3.1.1}By Lemma
3.1.2 we may assume $d$ divides $r$.
We consider the canonical
projection $\gamma : \AA^{m}_{k} \setminus \{0\}
\rightarrow \PP^{m - 1}_{k}$ and denote
by $\bar D$ the image of $D \setminus \{0\}$ in $\PP^{m - 1}_{k}$.
Let $h : Y \rightarrow \PP^{m - 1}_{k}$ be a resolution of $\bar D$ (in the
sense of 2.2). As in 2.2 we denote by $E_{i}$, $i \in J$,
the irreducible (smooth) components
of
$(h^{-1} (D))_{\rm red}$. We define similarly integers $N_{i}$ and
$\nu_{i}$,
and $E_{I}$, $J_{d}$, $U_{d}$, etc.

We denote by $U^{j}$ the open $x_{j} \not= 0$ in $\PP^{m - 1}_{k}$. The
restriction of $h$ to $h^{-1} (U_{j})$ is a resolution of
$f_{j} = \frac{f}{x_{j}^{r}}$ in $U_{j}$.
For $Z$ locally closed in $U_{d}$,
$[(Z \cap h^{-1} (U_{j}))_{f_{j}, \alpha}]$ has been defined in 2.2, and
by Lemma 1.4.1
$$
[(Z \cap h^{-1} (U_{j}) \cap h^{-1} (U_{j'}))_{f_{j}, \alpha}]
=
[(Z \cap h^{-1} (U_{j}) \cap h^{-1} (U_{j'}))_{f_{j'}, \alpha}].
$$
Thus we may define
without ambiguity
$$
[Z_{f, \alpha}] = \sum_{j} [(Z \cap h^{-1} (U_{j}))_{f_{j}, \alpha}]
-
\sum_{j \not= j'} [(Z \cap h^{-1} (U_{j}) \cap h^{-1} (U_{j'}))_{f_{j}, \alpha}]
+ \cdots.
$$

Set $E_{I}^{(d)} = E_{I} \setminus \cup_{\{j\} \notin J_{d}} E_{j}$.

\begin{prop}Assume
$d$ divides $r$.
For any character $\alpha$ of $\mu_{d} (k)$ of order $d$,
$$
\int_{X} (f^{s}, \alpha) = 
\frac{(\LL - 1) \, \LL^{- m}}{1 - \LL^{-rs -m}}
\sum_{I \in J_{d}}
[(E_{I}^{(d)})_{f, \alpha}] \prod_{i \in I}
\Bigl(
\frac{(\LL - 1) \, \LL^{-N_{i} s - n_{i}}}
{1 - \LL^{-N_{i} s - n_{i}}}
-1\Bigr).
$$
\end{prop}

\noindent{\em Proof.} Let us first remark that $$\int_{X} (f^{s}, \alpha) =
\frac{1}{1 - \LL^{-rs -m}} \int_{X \setminus \{0\}} (f^{s}, \alpha).$$
Indeed, by homogeneity of $f$, multiplication by $t$ induces an
isomorphism
between $Z_{n, f, \AA^{m}_{k}}$ and
$Z_{n + r, f, \{0\}}$ (notations of 2.1),
from which one deduces the relation
$$\LL^{m}   [X_{n + r, f, \{0\}},  \bar f^{\ast}\cL_{\alpha}]
=
\LL^{rm}  
[X_{n, f, \AA^{m}_{k}}, 
\bar f^{\ast}\cL_{\alpha}],$$ and the equality follows.

Write $X \setminus \{0\}$ as the disjoint union of the $W^{j}$'s, for
$1 \leq j \leq m$, with
$W^{j} = \{ x \in \AA^{m}_{k} \, \bigr| \,
x_{i} = 0 \, \,  \hbox{for} \, \,
i < j \, \, \hbox{and} \, \, x_{j} \not=0 \}.$
Now $\gamma (W^{j}) \subset U^{j}$ and
the restriction of $\gamma$ to $W^{j}$ is a trivial
fibration onto its image,
with fibre $\GG_{m, k}$. As the valuation of $f (\varphi (t))$ only
depends on $\gamma (\varphi (t))$
we deduce
$$
\int_{W^{j}} (f^{s}, \alpha) = (1 -\LL^{-1}) \int_{\gamma (W^{j})} (f_{j}^{s}, \alpha).
$$
Since $Y$ is the disjoint union of the
$h^{-1}(\gamma (W^{j}))$'s, we deduce from Theorem 2.2.1, by adding up,
that
$$
\int_{X \setminus \{0\}} (f^{s}, \alpha)
=
(1 -\LL^{-1}) \, \LL^{- (m-1)}
\sum_{I \in J_{d}}
[(E_{I}^{\circ})_{f, \alpha}] \prod_{i \in I}
\frac{(\LL - 1) \, \LL^{-N_{i} s - n_{i}}}
{1 - \LL^{-N_{i} s - n_{i}}}.
$$
The result follows, because
\begin{multline*}
\sum_{I \in J_{d}}
[(E_{I}^{\circ})_{f, \alpha}] \prod_{i \in I}
\frac{(\LL - 1) \, \LL^{-N_{i} s - n_{i}}}
{1 - \LL^{-N_{i} s - n_{i}}}
=\\
\sum_{I \in J_{d}}
[(E_{I}^{(d)})_{f, \alpha}] \prod_{i \in I}
\Bigl(\frac{(\LL - 1) \, \LL^{-N_{i} s - n_{i}}}
{1 - \LL^{-N_{i} s - n_{i}}}
-1\Bigr).\qed
\end{multline*}

\medskip By Proposition 3.2.1 we may write
$$\int_{X} (f^{s}, \alpha) = A \sum_{I \in J_{d}} [(E_{I}^{(d)})_{f,
\alpha}] \prod_{i \in I} B_{i},$$
with $$A = \frac{(\LL - 1) \, \LL^{- m}}
{1 - \LL^{-rs -m}}
\quad \hbox{and} \quad B_{i} = 
\frac{(\LL - 1) \, \LL^{-N_{i} s - n_{i}}}
{1 - \LL^{-N_{i} s - n_{i}}}
-1.$$
Remark that $A^{\vee} = \LL^{-rs} \LL^{m - 1} A$ and
$B_{i}^{\vee} = \LL^{-1} B_{i}$.
When $\alpha = 1$ the result follows because
$E_I$ being proper and smooth $\chi_c (E_I)^{\vee}
= \LL^{\vert I \vert - m + 1}
\chi_c (E_I)$. For general $\alpha$ the result
follows from the following lemma.\hfill$\qed$

\begin{lem}Assume 1.5.2 holds. For any $I$ in $J_{d}$,
$$[(E_{I}^{(d)})_{f, \alpha}]^{\vee}
=
\LL^{\vert I \vert - m + 1}
[(E_{I}^{(d)})_{f, \alpha^{-1}}].
$$
\end{lem}

\begin{proof}Let $Z$ be locally closed in $U_d$. If on a neighborhood of
$Z$ we may write
$f \circ h = u v^d$ with $u$ non vanishing, 
$[Z, u^{\ast} \cL_{\alpha}] = \chi_c ({\widetilde Z_{u, d}}, \alpha)$,
where $\widetilde Z_{u, d}$ is the cyclic cover defined in 1.4.
In general the covers $\widetilde Z_{u, d}$
can be glued together (cf. the proof of Lemma 1.4.1) to give a Galois cover
$\widetilde Z_d \rightarrow Z$ with group
$\mu_d (k)$, such that 
$[Z_{f, \alpha}] = \chi_c ({\widetilde Z_{d}}, \alpha)$,
for any character $\alpha$ of order $d$.

Set $W^{\circ} = E_I^{(d)}$ and $W = E_I$.

\begin{lem}The cover $\widetilde W_d^{\circ}
\rightarrow W^{\circ}$ extends to a ramified
$\mu_d (k)$-cover $\pi : 
\widetilde W_d
\rightarrow W$ which satisfies the following conditions.
\begin{enumerate}
\item [(1)]The scheme
$\widetilde W_d$ is proper and is locally for the Zariski topology
quotient of a smooth scheme $X$ by a finite abelian group $G$,
the $\mu_d (k)$-action on $\widetilde W_d$ being 
induced from a $\mu_d (k)$-action on $X$ commuting with the $G$-action.
\item[(2)]The morphism $\pi$ ramifies on
$\widetilde W_d \setminus \widetilde W_d^{\circ}$
and the $\mu_d (k)$-action on
the restriction of $\pi$
to $\widetilde W_d \setminus \widetilde W_d^{\circ}$
factors locally for the Zariski topology
trough a
$\mu_{d'} (k)$-action, for some $d' < d$ dividing
$d$.
\end{enumerate}
\end{lem}

\begin{proof}Let $x$ be a closed point of $W \setminus W^{\circ}$.
On a Zariski neighborhood $\Omega'$ of $x$ in $Y$,
$$f = u \prod_{i \in J_x} h_i^{N_i} v^d$$
with $u$ nonvanishing on $\Omega'$, 
$J_x = \{i \in J \setminus I \bigm \vert x \in E_i, d \not\vert N_{i}\}$ and $h_i$
local equations for $E_i$ near $x$. Put $\Omega = \Omega' \cap W$.
Since $\widetilde W^{\circ}_{d \vert \Omega \cap W^{\circ}}$ is
given by $y^d = u \prod_{i \in J_x} h_i^{N_i}$ in
$\GG_m \times (\Omega \cap W^{\circ})$,
we may extend $\widetilde W^{\circ}_{d \vert \Omega \cap W^{\circ}}
\rightarrow \Omega \cap W^{\circ}$ to
$\widetilde W_{d \vert \Omega}
\rightarrow \Omega$ by taking
$\widetilde W_{d \vert \Omega}$ to be
the normalization of the subscheme
$\bar W_{d \vert \Omega}$ of
$\AA^1 \times \Omega$ given
by $y^d = u \prod_{i \in J_x} h_i^{N_i}$, and
the $\mu_d (k)$-action
extends naturally.
The schemes $\widetilde W_{d \vert \Omega}
\rightarrow \Omega$ glue together
(cf. the proof of Lemma 1.4.1)
to give a scheme $\widetilde W_{d}
\rightarrow W$ with $\mu_d (k)$-action.

Let $d'$ be the gcd of $d$ and the $N_i$'s, $i \in J_x$.
We have $d' < d$. Locally for the \'{e}tale topology near
$x$, $\widetilde W_d$ is the disjoint union of the normalizations
of
$y^{d / d'} = \varepsilon^j u' \prod_{i \in J_x} h_i^{N_i / d'}$,
for $1 \leq j \leq d'$, for $\varepsilon$ a fixed
primitive $d'$-th root of unity and $u'$ such that
$u'{}^{d'} = u$. This implies that on a Zariski neighborhood
of $x$ in $W \setminus W^{\circ}$, the $\mu_d (k)$-action on
$\pi_{\vert \widetilde W_d
\setminus \widetilde W_d^{\circ}} : 
\widetilde W_d
\setminus \widetilde W_d^{\circ}
\rightarrow W \setminus W^{\circ}$ factors through a
$\mu_{d'} (k)$-action, since the normalization of a local complete
domain is again local.
We still have to verify that
locally for the Zariski topology $\widetilde W_d$ is the
quotient of a smooth scheme $X$ by a finite abelian group $G$,
the $\mu_d (k)$-action being 
induced from a $\mu_d (k)$-action on $X$
commuting with the $G$-action.
It is enough to check this for the scheme
$\widetilde W_{d \vert \Omega}$
which is
the normalization of the subscheme
$\bar W_{d \vert \Omega}$ of
$\AA^1 \times \Omega$ given
by $y^d = u \prod_{i \in J_x} h_i^{N_i}$.
We may assume $u = 1$. Indeed, consider the \'{e}tale cyclic
cover of degree $d$,
$p : \Omega' \rightarrow \Omega$ given by $u = u'{}^{d}$.
Since $\widetilde W_{d \vert \Omega}$ is 
the quotient of
the normalization 
of the subscheme
of
$\AA^1 \times \Omega'$ given
by $y^d = u'{}^{d} \prod_{i \in J_x} (h_i \circ p)^{N_i}$,
we are done by using the isomorphism
$\AA^1 \times \Omega' \rightarrow \AA^1 \times \Omega'$
given by $(y, x) \mapsto (y u'{}^{- 1}, x)$.
When $u = 1$
the scheme $\widetilde W_{d \vert \Omega}$ 
is the disjoint union of the normalizations
of
$y^{d / d'}=\varepsilon^j \prod_{i \in J_x} h_i^{N_i / d'}$,
for $1 \leq j \leq d'$, for $\varepsilon$ a fixed
primitive $d'$-th root of unity
and $d'$ the gcd of $d$ and the $N_i$'s, $i \in J_x$.
Hence we may finally assume that 
$\widetilde W_{d \vert \Omega}$ 
is the normalization 
of the scheme $\bar W_{d \vert \Omega}$ given by
$y^{d} = \prod_{i \in J_x} h_i^{N_i}$ in
$\AA^1 \times \Omega$,
and that the gcd of $d$ and the $N_i$'s, $i \in J_x$, is 1.
Now consider the subscheme $W'$ of
$\AA^1 \times \Omega \times \AA^{\vert J_x \vert}$
given by
$h_i = t_i^d$ and $y = \prod_{i \in J_x} t_i^{N_i}$.
It is easily seen that $W'$ is smooth. Let us denote by $\pi :
W' \rightarrow \bar W_{d \vert \Omega}$ the morphism
given by
$(y, x, (t_i)) \mapsto (y, x)$ and by $G$ the 
kernel of the morphism $\mu_d (k)^{\vert J_x \vert} \rightarrow
\mu_d (k)$ given by $(\xi_i) \mapsto \prod_{i \in J_x} \xi_i^{N_i}$.
The canonical action of $G$ on $\AA^{\vert J_x \vert}$ induces
an action on $W'$ for which the morphism $\pi$
is equivariant, hence $\pi$ factorizes through a morphism
$\pi_G : W' / G \rightarrow \bar W_{d \vert \Omega}$. Since
$\pi_G$ is of degree 1, the result follows.\end{proof}

\medskip We are now able to finish the proof. By Lemma 3.2.3 (2) and
Proposition 1.3.3 (2),
we have
$$[(E_{I}^{(d)})_{f, \alpha}] = \chi_{c} (\widetilde W_{d}, \alpha)
\quad \hbox{and} \quad [(E_{I}^{(d)})_{f, \alpha^{-1}}] = \chi_{c} (\widetilde W_{d},
\alpha^{-1}),$$
hence the result follows from the following proposition.\end{proof}

\begin{prop}Assume 1.5.2 holds.
Let $W$ be a proper $G$-scheme of pure dimension $m$, 
with $G$ a finite abelian group.
Assume that, locally for the Zariski topology,
$W$ is isomorphic as a $G$-scheme to a quotient
$X / H$ with $H$ a finite abelian group and $X$ a smooth 
$G \times H$-scheme.
Then, for any character $\alpha$ of 
$G$,
$$
\chi_{c} (W, \alpha)^{\vee}
=
\LL^{- m} 
\chi_{c} (W, \alpha^{-1}).
$$
\end{prop}

\begin{proof}By additivity of Euler characteristics (Theorem 1.3.1 (3) and
1.3.2 (6))
and by
Theorem 1.3.2 (5), we are reduced to prove that if $W$
is a  $G$-scheme of pure dimension $m$
which
is isomorphic as a $G$-scheme to a quotient
$X / H$ with $H$ a finite abelian group and $X$ a smooth 
$G \times H$-scheme, then, for any character $\alpha$ of 
$G$,
$
\chi (W, \alpha)^{\vee}
=
\LL^{- m} 
\chi_{c} (W, \alpha^{-1})$.
This follows directly
from Assertion 1.5.2 and Theorem 1.3.2 (4).\end{proof}

\medskip
Let us denote by $\LL'$ the image of $\LL$ by the morphism
$$\varphi : K_0 (\cM_{k, \QQ[\zeta_{d}]}) \longrightarrow K_0 (DM_{gm}
(k)_{\QQ[\zeta_{d}]}).$$ One then defines
a ring
$K_0 (DM_{gm}
(k)_{\QQ[\zeta_{d}]}) \, [\LL'{}^{s}, \LL'{}^{-s}]_{\rm loc}$
similarly as we defined the ring
$K_{0} (\cM_{k, \QQ[\zeta_{d}]}) \, [\LL^{s}, \LL^{-s}]_{\rm loc}$,
and $\varphi$ extends to a morphism
$K_{0} (\cM_{k, \QQ[\zeta_{d}]}) \, [\LL^{s}, \LL^{-s}]_{\rm loc}
\rightarrow
K_0 (DM_{gm}
(k)_{\QQ[\zeta_{d}]}) \, [\LL'{}^{s}, \LL'{}^{-s}]_{\rm loc}$,
which
we still denote by $\varphi$.

\begin{theorem}
For any character $\alpha$ of $\mu_{d} (k)$ of order $d$,
$$\varphi \Bigl(
\Bigl(\int_{X} (f^{s}, \alpha)\Bigr)^{\vee}\Bigr)
= \varphi \Bigl(
\LL^{- rs}
\int_{X} (f^{s}, \alpha^{-1})\Bigr)
$$
in
$K_0 (DM_{gm}
(k)_{\QQ[\zeta_{d}]}) \, [\LL'{}^{s}, \LL'{}^{-s}]_{\rm loc}$.
\end{theorem}

\begin{proof}The proof is the same as the one of Theorem 3.2.1,
using Corollary 1.5.4 instead of Assertion 1.5.2.\end{proof}

\section{Limit for $s \rightarrow  - \infty$ and nearby cycles}
\subsection{}We consider in this section the subring
$K_{0} (\cM_{k, \QQ[\zeta_{d}]}) \, [\LL^{-s}]_{\rm loc}''$
of
$K_{0} (\cM_{k, \QQ[\zeta_{d}]}) \, [[\LL^{-s}]]_{\rm loc}$
generated by the subring 
$K_{0} (\cM_{k, \QQ[\zeta_{d}]})$
and the series
$\LL^{-Ns - n} \, (1 - \LL^{-Ns - n})^{-1}$,
for $N$ and $n$ in $\NN \setminus \{0\}$.

\begin{lem}There is a well defined ring homomorphism
$${\rm CT} :
K_{0} (\cM_{k, \QQ[\zeta_{d}]}) \, [\LL^{-s}]_{\rm loc}''
\longrightarrow
K_{0} (\cM_{k, \QQ[\zeta_{d}]})$$ which induces the identity
on
$K_{0} (\cM_{k, \QQ[\zeta_{d}]})$, and which sends the series
$\LL^{-Ns - n} \, (1 - \LL^{-Ns - n})^{-1}$ to $-1$.
\end{lem}

\begin{proof}Similar to the one for the constant term
of the power series expansion in $T^{-1}$ of usual 
rational functions of degree $\leq 0$ in $T$, considering $\LL^{- s}$ as
a
variable $T$. (Heuristically this amounts to taking
the ``limit'' for $s \rightarrow - \infty$.)\end{proof}

\medskip
\begin{definition}Let $X$ be a smooth and
connected $k$-scheme of finite type of
dimension $m$, $f : X \rightarrow
\AA^{1}_{k}$ be a morphism and $x$ be a closed
point of $f^{-1} (0)$.
Let $\alpha$ be a character
of $\mu_{d} (k)$ of order $d$. 
We set 
$$
S_{\alpha, x} := \frac{\LL^{m}}{1 - \LL}
{\rm CT} \int_{\{x\}} (f^{s}, \alpha).
$$
By Theorem 2.2.1 and Lemma 4.1.1, 
$S_{\alpha, x}$ is an element of 
$K_{0} (\cM_{k, \QQ[\zeta_{d}]})$ which is well defined modulo
$(\LL - 1)$-torsion, cf. remark 1.2.3.
Furthermore, for any resolution of $f$, 
$$
S_{\alpha, x} =
\sum_{I \in J_{d}}
[(E_{I}^{\circ} \cap h^{-1} (x))_{f, \alpha}] 
(1 - \LL)^{\vert I \vert -1},
$$
modulo
$(\LL - 1)$-torsion.
\end{definition}

Remark that, for almost all $d$, $S_{\alpha, x} = 0$. 
We now assume for simplicity
that $k$ contains all roots of unity and that the group of 
roots of unity in $k$
is
embedded in $\CC$.
Hence all the groups $\widehat {\mu_{d} (k)}$ are canonically 
embedded in
$\QQ / \ZZ$. We denote by $\gamma$ the section $\QQ / \ZZ \rightarrow
[0, 1)$ and by $i_{x}$ the inclusion of $\{x\}$ in
$f^{-1} (0)$.

\medskip

We believe that $S_{\alpha, x}$ is the ``motivic incarnation''
of $\chi_{c} (i_{x}^{\ast} R \psi_{f, \alpha})$.
Here
$R \psi_{f, \alpha}$ denotes the eigenspace of nearby cycles
for the eigenvalue $\exp (2 \pi i
\gamma (\alpha))$ of the semi-simple part of the monodromy.
We will verify in the next subsection that this is true for the $\CC$-Hodge
realization.

\subsection{Hodge realization}We will
use freely the theory of mixed Hodge modules
developped by M.Saito in \cite{Sa1}, \cite{Sa3}. In particular, for $X$ a
scheme of finite type over $\CC$, we denote by
${\rm MHM} (X)$ the abelian category of mixed Hodge
modules on $X$. 
In the definition of mixed Hodge modules it is required that the
underlying
perverse sheaf is defined over $\QQ$. To allow  some more
flexibility
we will also use the category ${\rm MHM'} (X)$ of bifiltred
$\cD$-modules on $X$ which are direct factors of objects of
${\rm MHM} (X)$ as bifiltred
$\cD$-modules. We denote
by $D^{b} ({\rm MHM} (X))$ and $D^{b} ({\rm MHM'} (X))$
the corresponding derived
categories.

Let $f : X \rightarrow \AA^{1}_{\CC}$ be a morphism.
We denote by $\psi^{H}_{f}$ and $\phi^{H}_{f}$
the nearby and vanishing
cycle functors for mixed Hodge modules as defined in \cite{Sa3}
and $T_{s}$
the semi-simple part of the monodromy operator.
One should note that $\psi^{H}_{f}$ and $\phi^{H}_{f}$
on mixed Hodge modules correspond to $\psi_{f}[-1]$ and $\phi_{f}[-1]$
on the underlying perverse sheaves.
If $M$ is a mixed Hodge module on $X$ 
we denote by $\psi^{H}_{f, \alpha} M$ the object
of
${\rm MHM'} (X)$ which corresponds to
the eigenspace of $T_{s}$ for the
eigenvalue
$\exp (2 \pi i \gamma (\alpha))$.
These definitions extend to the Gro\-then\-dieck group of the
abelian category ${\rm MHM}' (X)$.

Let us recall the definition of {\it complex}
mixed Hodge structures.
A $\CC$-Hodge structure of weight $n$ is just a finite
dimensional bigraded vector space
$V = \bigoplus_{p + q = n} V^{p, q}$, or, equivalently,
a finite
dimensional vector space $V$ with decreasing filtrations $F^{\cdot}$ and
${\overline F}^{\cdot}$ such that $V = F^{p} \oplus {\overline F}^{q}$
when $p + q = n + 1$. A mixed
$\CC$-Hodge structure is
a finite
dimensional vector space $V$ with an increasing filtration $W$
and
decreasing filtrations $F^{\cdot}$ and
${\overline F}^{\cdot}$ 
which induce for each $n \in \ZZ$
a $\CC$-Hodge structure of weight $n$
on ${\rm Gr}^{W}_{n}V$.
We denote by
$K_{0} ({\rm MHS}_{\CC})$ the Grothendieck group of the
abelian category of complex
mixed Hodge
structures. 
The Hodge realization functor induces a morphism
$H : K_{0} (\cM_{\CC, \CC}) \rightarrow K_{0} ({\rm MHS}_{\CC})$.
Remark that $H (\LL) = \CC (-1)$ and that $H$ kills
$(\LL - 1)$-torsion (cf. remark 1.2.3).

For any object $K$ of $D^{b} ({\rm MHM} (X))$
we denote by $\chi_{c} (X, K)$ the class of
$Rp_{!} (K)$ in $K_{0} ({\rm MHS}_{\CC})$, where $p$ is the projection onto
${\rm Spec} \, \CC$. Clearly this definition may be extended to
$D^{b} ({\rm MHM'} (X))$.

If $X$ is smooth and connected of dimension $m$, we denote by 
$\CC_{X}^{H} [m]$ the trivial variation of Hodge structure of
weight 0.

\begin{theorem}Let $X$ be a smooth and connected $\CC$-scheme of finite type of
dimension $m$, $f : X \rightarrow
\AA^{1}_{\CC}$ be a morphism and $x$ be a closed
point of $f^{-1} (0)$.
The following equality holds
$$
H (S_{\alpha, x}) = (- 1)^{m - 1} \chi_{c} (i_{x}^{\ast} \psi_{f,
\alpha}^{H} \CC_{X}^{H} [m]).
$$
\end{theorem}

\begin{proof}We will use a resolution $h : Y \rightarrow X$
of $f$. We set $g = f \circ h$, $Y_{0} = g^{-1} (0)_{\rm red}$,
$Y^{\ast} = Y \setminus Y_{0}$, $g' = g_{\vert Y^{\ast}}$ and we denote
by
$j : Y^{\ast} \hookrightarrow Y$ the inclusion morphism.
We also set $D = h^{-1} (x)_{\rm red}$,
$D_{I} = D \cap E_{I}$, $D_{I}^{\circ} = D \cap E_{I}^{\circ}$,
and denote
by
$i_{D} : D \hookrightarrow Y_{0}$ the inclusion morphism.
In the derived category
$D^{b} ({\rm MHM'} (\GG_{m, \CC}))$ we have a decomposition
$$
R[d]_{\ast} \CC^{H} \simeq \oplus_{0 \leq j < d} \, \cK_{\frac{j}{d}},
$$
where the underlying sheaf of $\cK_{\frac{j}{d}}$
has monodromy
$\exp (2 \pi i \frac{j}{d})$ at the origin.
We set $\cF_{\alpha} := j_{\ast} g'{}^{\ast} \cK_{- \gamma (\alpha)}$.

\begin{lem}\begin{enumerate}
\item[(1)]The sheaf underlying $\cF_{\alpha}$ is locally constant of
rank 1 on $U_{d}$.
\item[(2)]For any $Z$ locally closed in $U_{d}$,
$$
\chi_{c} (Z, \cF_{\alpha}) = H ([Z_{f, \alpha}]).
$$
\item[(3)]For $i \geq 0$, $R^{i}j_{\ast} g'{}^{\ast} \cK_{- \gamma
(\alpha)}$
is zero outside $U_{d}$.
\item[(4)]For $i \geq 0$, $R^{i}j_{I \ast} \cF_{\alpha \vert E_{I}^{\circ}}$
is zero outside $E_{I} \cap U_{d}$,
with $j_{I} : E_{I}^{\circ} \hookrightarrow E_{I}$ the inclusion
morphism. 
\end{enumerate}
\end{lem}

\begin{proof}Assertion (1) follows from \cite{D2} Proposition 3.1.
Assertion (2) follows from the fact that, 
if locally $g = u v^{d}$ with $u$ non
vanishing on $Z$, then $u^{\ast} \cK_{- \gamma (\alpha)}$ is isomorphic
to
$\cF_{\alpha \vert Z}$. The proof of (3) and (4) are completely similar
to the one of
\cite{D2} Lemma 3.2.\end{proof}

The functor $R^{j} (h_{\vert Y_{0}})_{\ast} :
D^{b} ({\rm MHM} (Y_{0})) \rightarrow {\rm MHM}
(f^{-1} (0)_{\rm red})$ being a cohomological functor, we have a
spectral sequence (cf. \cite{Sa4} (2.14.3))
$$
E_{2}^{pq} = R^{p} (h_{\vert D})_{\ast} H^{q}
i_{D}^{\ast} (\psi^{H}_{g}  \CC^{H}_{Y} [m]) \Longrightarrow
H^{p + q}i_{x}^{\ast} (\psi^{H}_{f}  \CC^{H}_{X} [m]),
$$
which is $T_{s}$-equivariant. Hence, because the restriction of $h$ to
$D$ is proper, it is enough, by Proposition 4.2.2 (2),
to prove the following proposition.\end{proof}

\begin{prop}With the previous notations, the following equality holds
$$
\chi_{c} (D, i_{D}^{\ast} (\psi^{H}_{g, \alpha}  \CC^{H}_{Y} [m]))
=
(-1)^{m - 1} \sum_{I \in J_{d}} \chi_{c} (D_{I}^{\circ}, \cF_{\alpha})
(1 - \CC (-1))^{\vert I \vert -1}.
$$
\end{prop}

\noindent{\em Proof.} Let $N$ be the logarithm of the unipotent part of the monodromy
and let $P_{N}$ denote the primitive part with respect to $N$.
We have the primitive decomposition
\numberwithin{equation}{theorem}
\begin{multline}
\bigoplus_{j}
{\rm Gr}_{j}^{W}
\psi^{H}_{g, \alpha}  \CC^{H}_{Y} [m]
\simeq\\
\bigoplus_{k \geq 0} \bigoplus_{i = 0}^{k}
\Bigl[N^{i} P_{N} {\rm Gr}_{m - 1 + k}^{W}
\psi^{H}_{g, \alpha}  \CC^{H}_{Y} [m]\Bigr] (i).
\end{multline}

The proposition will follow from the following lemma.

\begin{lem}For $k \geq 0$, there is a canonical isomorphism
$$P_{N} {\rm Gr}_{m - 1 + k}^{W}
\psi^{H}_{g, \alpha}  \CC^{H}_{Y} [m]
\simeq
\bigoplus_{I \in J_{d}, \, \vert I \vert = k + 1} {\rm IC}_{E_{I}}
\cF_{\alpha \vert E_{I}^{\circ}} (-k),
$$
with ${\rm IC}_{E_{I}}
\cF_{\alpha \vert E_{I}^{\circ}}$ the intersection cohomology module on $E_{I}$ with
coefficients
in $\cF_{\alpha \vert E_{I}^{\circ}}$.
\end{lem}

\begin{proof}The statement for the underlyings perverse sheaves
is essentially Lemma 2.13 of \cite{Sa4}. The preprint \cite{Sa4} being
unpublished, we reproduce the argument for the convenience of the
reader. Let us denote by $S$ the functor which to an object of
${\rm MHM'}$ associates its underlying perverse sheaf.
By \cite{Sa2} \S\kern .15em 3 we have a canonical isomorphism
\begin{multline}
S \Bigl[\psi^{H}_{g, \alpha}  \CC^{H}_{Y} [m]\Bigr]
\simeq\\
\limind {\rm Ker} \, (j_{!} g'{}^{\ast} \cK_{- \gamma(\alpha), k} [m]
\longrightarrow
j_{\ast} g'{}^{\ast} \cK_{- \gamma (\alpha), k} [m]), 
\end{multline}
where, for $k \geq 0$, $\cK_{- \gamma (\alpha), k}$
is a local system of rank $k + 1$ on $\GG_{m,
\CC}$ whose monodromy at the origin has a unique Jordan block and
eigenvalue
$\exp (- 2 \pi i \gamma (\alpha))$, and 
$\cK_{- \gamma (\alpha), 0}$ is isomorphic to 
$S (\cK_{- \gamma (\alpha)})$.
We have a canonical isomorphism
$$
{\rm Gr}_{m - \vert I \vert}^{W}j_{!} g'{}^{\ast} \cK_{- \gamma (\alpha), k}
[m]_{\vert E_{I}^{\circ}}
\simeq
S \Bigl[\cF_{\alpha} [m - \vert I \vert]_{\vert E_{I}^{\circ}}\Bigr]
$$
from which one  deduces the canonical isomorphism
$$
{\rm Gr}_{m - k}^{W}j_{!} g'{}^{\ast} \cK_{- \gamma (\alpha), k}
[m]
\simeq
S \Bigl[\bigoplus_{I \in J_{d}, \, \vert I \vert = k} {\rm IC}_{E_{I}}
\cF_{\alpha \vert E_{I}^{\circ}}\Bigr],
$$
for $k \geq 0$, and 
one can check that the natural morphism 
$$j_{!} g'{}^{\ast} \cF_{\alpha} [m] \longrightarrow \psi^{H}_{g, \alpha}
\CC^{H}_{Y} [m]$$
deduced from (4.2.4.1)
induces an isomorphism
$$
{\rm Gr}_{m - 1 -k}^{W}j_{!} g'{}^{\ast} \cK_{- \gamma (\alpha), k}
[m]
\simeq
S \Bigl[N^{k} P_{N} {\rm Gr}_{m - 1 + k}^{W}
\psi^{H}_{g, \alpha}  \CC^{H}_{Y} [m]\Bigr]
$$ by loc.cit. Hence we have a canonical isomorphism
$$S \Bigl[P_{N} {\rm Gr}_{m - 1 + k}^{W}
\psi^{H}_{g, \alpha}  \CC^{H}_{Y} [m]\Bigr]
\simeq
S \Bigl[\bigoplus_{I \in J_{d}, \, \vert I \vert = k + 1} {\rm IC}_{E_{I}}
\cF_{\alpha \vert E_{I}^{\circ}} (-k)\Bigr].$$
The fact that this isomorphism lifts canonically to an isomorphism
between the corresponding objects
of
${\rm MHM'}$ follows from \cite{Sa1} 3.6.10 and 5.2.16 (cf. \cite{Sa1} p.990).\end{proof}

\medskip

By primitive decomposition (4.2.3.1),
\begin{equation*}
\begin{split}
\chi_{c} (D, i_{D}^{\ast} (\psi^{H}_{g, \alpha}  \CC^{H}_{Y} [m]))
&=
\sum_{k} \sum_{i = 0}^{k}
\chi_{c} (D, (i_{D}^{\ast} N^{i} P_{N} {\rm Gr}_{m - 1 + k}^{W}
\psi^{H}_{g, \alpha}  \CC^{H}_{Y} [m]) (i))\\
&=
\sum_{k}
\chi_{c} (D, i_{D}^{\ast}  P_{N} {\rm Gr}_{m - 1 + k}^{W}
\psi^{H}_{g, \alpha}  \CC^{H}_{Y} [m]) \, (\sum_{i = 0}^{k} \CC(i)).
\end{split}
\end{equation*}

By Lemma 4.2.4 we deduce
\begin{multline*}
\chi_{c} (D, i_{D}^{\ast} (\psi^{H}_{g, \alpha}  \CC^{H}_{Y} [m]))
=\\
\sum_{k} 
\sum_{I \in J_{d}, \vert I \vert = k + 1} \chi_{c} (D,
i_{D}^{\ast}
{\rm IC}_{E_{I}} \cF_{\alpha \vert E_{I}^{\circ}}) \, (\sum_{i = 0}^{k} \CC(-i)).
\end{multline*}

By Lemma 4.2.2  we have
\begin{equation*}
\begin{split}
\chi_{c} (D, i_{D}^{\ast} {\rm IC}_{E_{I}} \cF_{\alpha \vert E_{I}^{\circ}})
&=
\chi_{c} (D_{I}, j_{I \ast} \cF_{\alpha \vert E_{I}^{\circ}} [m - \vert
I \vert])\\
&=
(-1)^{m - \vert I \vert}
\sum_{{I' \supset I} \atop I' \in J_{d}} \chi_{c} (D_{I'}^{\circ},
\cF_{\alpha \vert D_{I'}^{\circ}}),
\end{split}
\end{equation*}
where $j_{I}$ denotes the inclusion $E_{I}^{\circ} \hookrightarrow E_{I}$.

Hence we deduce
$$
\chi_{c} (D, i_{D}^{\ast} (\psi^{H}_{g, \alpha}  \CC^{H}_{Y} [m]))
=
\sum_{I \in J_{d}}
\alpha_{I} \,
\chi_{c} (D_{I}^{\circ},
\cF_{\alpha \vert D_{I}^{\circ}}),
$$
with 
$$
\alpha_{I} 
=
\sum_{k = 0}^{\vert I \vert - 1}
(- 1)^{m - k - 1} {\vert I \vert \choose k + 1}
(1 + \CC(-1) + \cdots + \CC(-k)).
$$
The result follows, since, by the binomial theorem,
$$
\alpha_{I} 
=
(-1)^{m - 1} (1 - \CC (-1))^{\vert I \vert - 1}.\hfill\qed
$$

\theoremstyle{remark}
\newtheorem*{Rk}{Remark}\begin{Rk}Of course one deduces from Theorem 4.2.1
the analagous statement for the Betti realization, which
may also be proved directly using Lemma 4.2.2 and Leray's spectral
sequence for nearby cycles. See also \cite{D4} for a related result
concerning the \'{e}tale realization.\end{Rk}

\subsection{Relation with the Hodge spectrum}Let
us review the definition of Hodge spectrum
according to \cite{Sa5}. Let $H$ be a complex
mixed Hodge structure with an
automorphism
$T$ of order dividing $d$.
The Hodge spectrum of $(H, T)$ is defined as
${\rm HSp} (H, T) = \sum_{\alpha \in \frac{1}{d} \ZZ}  n_{\alpha}
t^{\alpha} \in
\ZZ [t^{- \frac{1}{d}}, t^{\frac{1}{d}}]$, with
$n_{\alpha} = {\rm dim \, Gr}^{p}_{F} H_{\lambda}$,
for $\lambda = \exp (2 \pi i \alpha)$ and $p = [\alpha]$,
where
$H_{\lambda}$ is the eigenspace of $T$ with eigenvalue $\lambda$,
and
$F$ is the Hodge filtration.
This definition extends to the Grothendieck group of the
abelian category of
complex mixed Hodge structures with an
automorphism
$T$ of order dividing $d$.
Remark that ${\rm HSp} (H (k), T) = t^{-k} {\rm HSp} (H, T)$, where
$(k)$ is the Tate twist.

We assume now $X$ to be a
smooth and connected scheme of finite type over $\CC$ and dimension $m$.
We set
$${\rm HSp} (f, x) = \sum_{j \in \ZZ} (-1)^{j}  {\rm HSp}
(H^{j}i_{x}^{\ast}
\psi^{H}_{f} \CC_{X}^{H} [m], T_{s})
$$
and
$${\rm HSp}' (f, x) = \sum_{j \in \ZZ} (-1)^{j}  {\rm HSp}
(H^{j}i_{x}^{\ast}
\phi^{H}_{f} \CC_{X}^{H} [m], T_{s}).
$$
We have
${\rm HSp} (f, x)
=
(-1)^{m - 1} + {\rm HSp}' (f, x)$. The relation with the spectrum
${\rm Sp} (f, x)$ as defined in \cite{St1} and \cite{Sa5} (which differs from
that of \cite{St2} by multiplication by $t$)
is 
$$
{\rm Sp} (f, x) = t^{m} \iota ({\rm HSp}' (f, x))
$$
where $\iota$ is the isomorphism such that $\iota (t^{\frac{1}{d}})
= t^{- \frac{1}{d}}$.

The following statement is a direct corollary of Theorem 4.2.1.

\begin{cor}Let $X$ be a smooth and connected $\CC$-scheme of finite type of
dimension $m$, $f : X \rightarrow
\AA^{1}_{\CC}$ be a morphism and $x$ be a closed
point of $f^{-1} (0)$.
The following equality holds
$$
\sum_{\alpha \in \QQ / \ZZ}  {\rm HSp} (H (S_{\alpha, x}), {\rm Id})
t^{\gamma (\alpha)}
= (- 1)^{m - 1} {\rm HSp} (f, x). \hfill \qed
$$
\end{cor}

\section*{Appendix}
Let $G$ be a finite group.
As in 1.3 we denote by ${\rm Sch}_{k, G}$ the category of  
separated schemes of finite type over $k$ with $G$-action
such that
the $G$-orbit of any closed point of $X$ is contained in
an affine open subscheme.
We also denote
by ${\rm Reg}_{k, G}$ the full subcategory of smooth schemes.
In this section we state some
analogs  with group action
of resolution statements used in \cite{G-N}.
In \cite{G-N} the  statements
were deduced from
Hironaka's theorems \cite{H1}. Here we  
will use instead Villamayor's
results in \cite{V2} \S\kern .15em 7
(which, as indicated in the introduction
of \cite{V2}, may also be deduced from
\cite{H2} together with \cite{V1}).

Let us begin by some elementary observations. By the very definition,
objects of ${\rm Sch}_{k, G}$ have covers
by $G$-stable affine subschemes. Furthermore,
such affine subschemes may be embedded equivariantly in
a smooth $G$-scheme. Indeed, if $X$ is an affine $G$-scheme,
it may be embedded in a smooth scheme $Z$, and the embedding
$X \hookrightarrow Z^G$ given by $x \mapsto (gx)_{g \in G}$
is equivariant. Similarly any object $X$ of
${\rm Sch}_{k, G}$ admits an equivariant compactification:
if $Z$ is any compactification of $X$, the closure of the image of
$X$
by
$x \mapsto (gx)_{g \in G}$ in $Z^G$ gives an equivariant
compactification.
Remark that the equivaraint Chow Lemma may be deduced
directly from the existence of equivariant compactifications by the
usual proof of the Chow Lemma.

By a proper relative isomorphism of $G$-schemes
$(\tilde X, \tilde Y) \rightarrow (X, Y)$, we mean the following
data~:
a proper $G$-morphism $f : \tilde X \rightarrow X$ between objects of
${\rm Sch}_{k, G}$, reduced closed $G$-stable
subschemes $\tilde Y$ and
$Y$ of $\tilde X$ and
$X$ respectively, such that $\tilde Y$ is the preimage of
$Y$ in $\tilde X$, and such that the restriction of $f$ to
$\tilde X \setminus \tilde Y$ is an $G$-isomorphism
onto $X \setminus Y$.
If moreover $X$ is smooth,
$Y$ is a closed subscheme which is
smooth and of smaller dimension,
and $f : \tilde X \rightarrow X$
is isomorphic to the  blowing up of $X$ along $Y$, we will say
$f$ is an elementary $G$-modification.

\theoremstyle{plain}
\newtheorem*{thA.1}{Theorem A.1}
\begin{thA.1}For any reduced object $X$ in
${\rm Sch}_{k, G}$,  there exists 
a proper relative isomorphism of $G$-schemes
$(\tilde X, \tilde Y) \rightarrow (X, Y)$
with 
$\tilde X$ smooth and ${\rm dim} \, Y$,
${\rm dim} \, \tilde Y < {\rm dim} \, X$.
\end{thA.1}

\begin{proof}Follows directly from \cite{V2} \S\kern .15em 7.\end{proof}

\medskip

We will also need the following equivariant
Chow-Hironaka Lemma.

\newtheorem*{Lemma A.2}{Lemma A.2}

\begin{Lemma A.2}Assume $G$ is a finite group.
Let $f : \tilde X \rightarrow X$ be a birational proper
morphism
in ${\rm Reg}_{k, G}$, inducing a 
birational proper
morphism on each irreducible component of $\tilde X$. Then there exists
a $G$-morphism $X' \rightarrow X$, which is the
composition of a finite sequence of elementary $G$-modifications
and
which factors through a proper birational $G$-morphism
$X' \rightarrow \tilde X$.
\end{Lemma A.2}

\begin{proof}By the Chow-Hironaka Lemma one can dominate $f$
by a projective birational map $X' \rightarrow X$ which is
the blowing up
of a sheaf of ideals
$\cJ$. Let $\cI = g_1 \cJ \cdot g_2 \cJ \cdots
g_n \cJ$, with $g_1, \ldots, g_n$ the elements of $G$.
Since the blowing up of $\cI$ dominates the one of $\cJ$,
it is enough to know that it is possible
to make $\cI$ principal
by a
composition of a finite sequence of elementary $G$-modifications,
a fact
which follows from \cite{V2} \S\kern .15em 7.\end{proof}

\medskip

We denote by ${\rm Reg}_{k, G}^2$ the
category of pairs $(X, U)$ with $X$ in ${\rm Reg}_{k, G}$
and $U$ an $G$-stable open subscheme of $X$ such
that the complement $D = X \setminus U$ is a $G$-stable 
divisor with normal crossings.

\newtheorem*{Lemma A.3}{Lemma A.3}

\begin{Lemma A.3}Assume $G$ is a finite group.
Let $U$ be object of ${\rm Reg}_{k, G}$
and $V$ be a smooth $G$-stable subscheme.
There exists a compactification
$X$ of $U$ such that
$(X, U)$ belongs to
${\rm Reg}_{k, G}^2$
and the closure of $V$ in $X$ has normal crossings with $X \setminus
U$.
\end{Lemma A.3}

\begin{proof}It is the same proof as the one in \cite{G-N} Lemme 3.5,
using \cite{V2} instead of
\cite{H1}, once remarked that $U$ admits a compactification 
with $G$-action.\end{proof}

\medskip

A morphism $f : (\tilde X, \tilde U) \rightarrow (X, U)$ in
${\rm Reg}_{k, G}^2$ will be called an elementary $G$-modification
if $f : \tilde X \rightarrow X$
is the blowing up
of $X$ along a smooth $G$-stable subscheme $Y$
having normal crossings with $X \setminus U$ and
furthermore $\tilde U =
f^{- 1} U$.
With this terminology, the analogue of Lemma 2.7 of \cite{G-N}
holds in the $G$-equivariant setting by using \cite{V2} instead of
\cite{H1}.

\medskip

\begin{proof}[Proof of Theorem 1.3.1 and 1.3.2]The construction
and proof of unicity of $\chi_c$ satisfying
conditions (1) and (2) in Theorem 1.3.1 is just done
the same
as for the analogue statement in \cite{G-N} (Corollaire 6.4), using the previous
equivariant versions of the resolution results used in \cite{G-N}.
Assertion (3) is easily deduced from (2).
Similarly, the construction of $\chi$
satisfying
conditions (1) to (4) in Theorem 1.3.2, and determined by the first
three, is done exactly as in \cite{G-N} (Corollaire 6.13). Assertion (5)
follows by construction, and (6) is easily deduced from (2) and
(5) by induction on dimension.\end{proof}

\bibliographystyle{amsplain}

\end{document}